\UseRawInputEncoding
\documentclass[a4paper,oneside]{article}
\usepackage{amsmath, amssymb, amsthm,graphicx}
\usepackage[font=small,labelfont=md,textfont=it]{caption}
\usepackage{floatrow,float}
\usepackage[titletoc, title]{appendix}
\usepackage[colorlinks,linkcolor=blue,citecolor=blue]{hyperref}
\usepackage{longtable}
\usepackage{pgfplots}
\usepackage{diagbox}
\usepackage{booktabs,makecell,multirow}
\usepackage[capitalise, nosort]{cleveref}
\usepackage{cases,color}
\usepackage{algorithm}
\usepackage{algorithmic}
\usepackage{amsmath,bm}
\usepackage{geometry}
\usepackage{extarrows}
\usepackage{tcolorbox}

\DeclareMathOperator{\D}{D}
\DeclareMathOperator{\I}{I}

\crefname{equation}{}{}
\crefname{lem}{Lemma}{Lemmas}
\crefname{thm}{Theorem}{Theorems}

\newcommand{\dual}[1]{\left\langle {#1} \right\rangle}

\newcommand{\jmp}[1]{{[\![ {#1} ]\!]}}
\newcommand{\nm}[1]{\left\lVert {#1} \right\rVert}
\newcommand{\snm}[1]{\left\lvert {#1} \right\rvert}
\newcommand{\ssnm}[1]
{
	\left\vert\kern-0.25ex
	\left\vert\kern-0.25ex
	\left\vert
	{#1}
	\right\vert\kern-0.25ex
	\right\vert\kern-0.25ex
	\right\vert
}

\makeatletter
\def\spher@harm#1{%
	\vbox{\hbox{%
			\offinterlineskip
			\valign{&\hb@xt@2\p@{\hss$##$\hss}\vskip.2ex\cr#1\crcr}%
		}\vskip-.36ex}%
}
\def\gshone{\spher@harm{.}}
\def\gshtwo{\spher@harm{.&.}}
\def\gshthree{\spher@harm{.&.&.}}
\let\gsh\spher@harm
\makeatother

\newtheorem{lemma}{Lemma}[section]

\newtheorem{remark}{Remark}[section]

\newtheorem{theorem}{Theorem}[section]

\newcounter{mnote}
\let\oldmarginpar\marginpar
\renewcommand\marginpar[1]
{\-\oldmarginpar[\raggedleft\footnotesize #1]{\raggedright\footnotesize #1}}

\makeatletter\def\@captype{table}\makeatother

\newfloatcommand{capbtabbox}{table}[][\FBwidth]
\begin{document}
	\title{
		\Large \bf Error estimation of a discontinuous Galerkin method for time fractional
		subdiffusion problems with nonsmooth data
		\thanks
		{
			This work was supported in part by National Natural Science Foundation
			of China (11771312).
		}
	}
	\author{
		Binjie Li \thanks{School of Mathematics, Sichuan University, Chengdu, 610064, China. Email: libinjie@scu.edu.cn},
		Hao Luo \thanks{School of Mathematical Sciences, Peking University, Beijing, 100871, China. Corresponding author. Email: luohao@math.pku.edu.cn},
		Xiaoping Xie \thanks{School of Mathematics, Sichuan University, Chengdu, 610064, China. Email: xpxie@scu.edu.cn} \\
	}
	\date{}
	\maketitle
	\begin{abstract}
	This paper is devoted to the numerical analysis of a piecewise constant discontinuous Galerkin method for time fractional subdiffusion problems. The  regularity of weak solution is firstly established by using variational approach and Mittag-Leffler function. Then several optimal error estimates are derived with low regularity data. Finally,
	numerical experiments are conducted to verify the theoretical results.
	\end{abstract}
	
	
	
\section{Introduction} \label{sec:1}

\setcounter{section}{1} \setcounter{equation}{0}

This paper considers the lowest-order discontinuous Galerkin (DG) method for the time fractional subdiffusion equation:
\begin{equation}
	\label{eq:model}
	\left\{
	\begin{aligned}
		u' - \D_{0+}^{1-\alpha} \Delta u & = f   &  & \text{in $ \Omega_T:=\Omega\times(0,T)$,}         \\
		u                                             & = 0   &  & \text{on $ \partial\Omega \times (0,T) $,} \\
		u(\cdot,0)                                    & = u_0 &  & \text{in $ \Omega $,}
	\end{aligned}
	\right.
\end{equation}
where $\alpha \in(0,1),\, T>0$ denotes the final time, $\Omega \subset \mathbb R^d $ ($ d=1,2,3$) is a convex polyhedral domain, $ f $ and $ u_0 $ are given data, and $ \D_{0+}^{1-\alpha} $ is a left-sided Riemann\textendash Liouville fractional differential operator (cf. Section~\ref{sec:2}).

Let us briefly review the efforts devoted to numerical treatments
of problem \eqref{eq:model}. By now there are mainly two types of methods, according to how
the time fractional derivative is approximated. The first type of schemes are based on finite difference formula, including L-type schemes \cite{Gao2011,Langlands2005,Li2019L1,Zhuang2008} and the Gr\"unwald--Letnikov (GL) scheme \cite{Mohebbi2013,Yuste2006,Yuste2005}. The L1 method has the accuracy $O(\tau^{1+\alpha})$ for $C^2$ solutions. Utilizing the superconvergence property at some particular points of the GL formula, Gao et al. \cite{Gao2015} constructed some finite difference schemes that achieve the accuracy $O(\tau^2)$ for $C^3$ solutions. We also mention that for $u_0\in\dot H^2(\Omega)$ and smooth $f$, quadrature rules with exponential rate $O(e^{-c\sqrt{N}})$ can be found in \cite{Thomee2010IMA,Thomee2010}, where $c>0$ is some constant and $N$ stands for the the number of quadrature points.

The second type of schemes adopt the time-stepping DG method \cite{Eriksson1985,schotzau_time_2000}
and are combined with graded temporal grids to conquer the singularity. In \cite{McLean2009Convergence},
McLean et al. has applied the piecewise constant DG method to problem \eqref{eq:model} and proved the error bound $O(\tau+\snm{\ln \tau} h^2)$ under the $L^\infty(0,T;L^2(\Omega))$-norm, with initial data $u_0\in\dot H^2(\Omega)$ and the following regularity assumption:
\begin{equation}\label{eq:regu-assum}
	\begin{aligned}
		\big\|u(t)\big\|_{\dot H^2(\Omega)} + t\big\|u'(t)\big\|_{\dot H^2(\Omega)}
		\leqslant {}&M &&\quad 0<t\leqslant T,\\
		t^{2-\alpha}\big\|u'(t)\big\|_{\dot H^2(\Omega)} + t^{3-\alpha}\big\|u''(t)\big\|_{\dot H^2(\Omega)}
		\leqslant {}&M t^{\sigma-1}&&\quad 0<t\leqslant T,
	\end{aligned}
\end{equation}
where $\sigma$ and $M$ are two positive constants. Besides, we list more works using the piecewise linear DG method: Mustapha et al. \cite{Mustapha2011Piecewise} proved that the temporal convergence order under the $L^\infty(0,T;L^2(\Omega))$-norm is $O(\tau^{1+\alpha})$; later in \cite{Mustapha2012Uniform}, the authors derived
the improved bound $O(\tau^{\min\{1.5+\alpha,\,2\}})$; in addition, they \cite{Mustapha2013} proved the rate $O(\snm{\ln \tau}\tau^{1+2\alpha})$, which yields superconvergence if $\alpha\in(1/2,1)$. We mention that the analyses in \cite{Mustapha2013,Mustapha2012Uniform} require stronger growth assumptions than \eqref{eq:regu-assum}. Mustapha \cite{Mustapha2015Time} also proposed an $hp$-version DG method for problem \eqref{eq:model} and established the suboptimal order $O(\tau^{\max\{k,2\}+(1-\alpha)/2})$,
where $k\geqslant 1$ is the polynomial degree.

It is worth to noticing an alternative form of \eqref{eq:model}:
\begin{equation}
	\label{eq:model-equi}
	\left\{
	\begin{aligned}
		\D_{0+}^{\alpha} (u-u_0)-\Delta u & = \widetilde{f}    &  & \text{in $ \Omega_T$,}         \\
		u                                             & = 0   &  & \text{on $ \partial\Omega \times (0,T) $,} \\
		u(\cdot,0)                                    & = u_0 &  & \text{in $ \Omega $,}
	\end{aligned}
	\right.
\end{equation}
where $\widetilde{f}= \I_{0+}^{1-\alpha}f $ with $ \I_{0+}^{1-\alpha} $ being a left-sided Riemann\textendash Liouville fractional integral operator (cf. Section~\ref{sec:2}). For $f=0$, both \eqref{eq:model} and \eqref{eq:model-equi} share the same solution that can be represented by the Mittag-Leffler function. For solution regularity and numerical analysis of problem \eqref{eq:model-equi}, especially for nonsmooth data, we refer the readers to \cite{ford_finite_2011,Li2018A,Li2019Numer,Li2020Numer,li_regu_2019,Li2018,mustapha_discontinuous_2014,stynes_2016,stynes_2017,wang2020,yan_analysis_2018,yang_time_2018}.

In a series of works \cite{McLean2010,Thomee2010IMA,Thomee2010}, by using Laplace transform, McLean et al. considered the regularity of problem \eqref{eq:model} and proved some growth estimates. To our best knowledge, no work has been proposed for investigating the weak solution to problem \eqref{eq:model} by variational approach. In this paper, for the case $u_0=0,\,f\neq 0$, we introduce a weak solution to problem \eqref{eq:model} via variational formulation, and prove that if $f\in L^2(0,T;\dot H^{-\beta}(\Omega))$ with $0\leqslant\beta \leqslant1$, then
\begin{equation*}
	\begin{aligned}
		{}&\nm{u}_{{}_0H^1(0,T;\dot H^{-\beta}(\Omega))} +
		\nm{u}_{{}_0H^{1-\alpha}(0,T;\dot H^{2-\beta}(\Omega))}
		\leqslant{} C_\alpha\nm{f}_{L^2(0,T;\dot H^{-\beta}(\Omega))}.
	\end{aligned}
\end{equation*}
For nonhomogeneous case: $f=0,\,u_0\neq0$, the weak solution is introduced and analyzed by Mittag-Leffler function.
Here we note that, instead of proving the growth estimate like \eqref{eq:regu-assum}, we show what kind of vector-valued Sobolev space the weak solution belongs to.

As  mentioned before, the error analyses of   most existing numerical methods require either smooth property or growth estimate of the true solution. It will be more challengable to establish error estimates with given low regularity data. Let us summarize few work that aims to fill in this gap. For a temporal semi-discretization with $ f= 0 $, McLean and
Mustapha \cite{McLean2015Time} derived that
\[
\nm{(u-u_\tau)(t_j)}_{L^2(\Omega)}
\lesssim
t_j^{-1} \tau \nm{u_0}_{L^2(\Omega)},
\]
by using the Laplace transform. For a spatial semi-discretization, by energy argument, Karaa et al.~\cite{Karaa2018} proved that, for $0<t\leqslant T$,
\[
\nm{(u\!-\!u_h)(t)}_{L^2(\Omega)}\!
\lesssim\!
h^2 t^{-\alpha(2-\delta)/2}\!
\left(
\nm{u_0}_{\dot H^\delta(\Omega)} \!+\!
\sum_{i=0}^2\int_0^Tt^i\nm{f^{(i)}(t)}_{\dot H^\delta(\Omega)} \!\mathrm{d}t
\right),
\]
where $ 0 \leqslant \delta \leqslant 2 $.

In this work, we shall derive several optimal error estimates for a piecewise constant DG method with nonsmooth data:
\begin{itemize}
	\item[\textbullet] if $ f = 0 $ and $ u_0 \in  L^2(\Omega)$, then
	\begin{equation}\label{eq:conv-u0-L2}
		\nm{u-U}_{L^2(\Omega_T)}
		\lesssim
		\big( \tau^{1/2} +h \big)
		\nm{u_0}_{L^{2}(\Omega)};
	\end{equation}
	\item[\textbullet] if $ u_0 = 0 $ and $ f \in L^2(\Omega_T) $, then
	\begin{align}
		(\tau^{1/2}\!+\!h)\!
		\nm{u\!-\!U}_{H^{\frac{1-\alpha}2}(0,T;\dot H^1(\Omega))}
		\!\!+\!\!
		\nm{u\!-\!U}_{L^2(\Omega_T)}
		\!\lesssim{}&\!
		\big( \tau \!+\!h^2 \big)\!
		\nm{f}_{L^2(\Omega_T)}.
		\label{eq:conv-f-L2-L2-intro}
	\end{align}
\end{itemize}
Note that \eqref{eq:conv-f-L2-L2-intro} is optimal  with respect to the solution regularity (cf. Remark~\ref{rem:opt-rate}) but \eqref{eq:conv-u0-L2} is optimal only for temporal discretization. Moreover, for the case $u_0=0$ with uniform temporal grid, by means of Laplace transform, we prove the following quasi-optimal results:
\begin{itemize}
	\item[\textbullet] if $ f \in L^2(\Omega_T) $, then
	\[
	\nm{u - U}_{L^\infty(0,T;L^{2}(\Omega))}
	\lesssim{}\snm{\ln \tau}
	\left(\tau^{1/2} +\epsilon_h
	h^{\min\{2,1/\alpha\}}\right)
	\nm{f}_{L^2(\Omega_T)},
	\]
	where $\epsilon_h = 1$ if $\alpha\neq1/2$ and $\epsilon_h = \sqrt{\snm{\ln h}}$ if $\alpha=1/2$;
	\item[\textbullet] if $ f \in {}_0H^{1/2}(0,T;L^2(\Omega)) $, then
	\[
	\nm{u - U}_{L^\infty(0,T;L^2(\Omega))}
	\lesssim{}\snm{\ln \tau}
	\left(\snm{\ln \tau}\tau+  h^2 \right)
	\nm{f}_{{}_0H^{1/2}(0,T;L^2(\Omega))}.
	\]
\end{itemize}

The rest of this paper is organized as follows. Section~\ref{sec:2} introduces some
Sobolev spaces and fractional calculus operators.
Section~\ref{sec:3} defines the weak solution to problem \eqref{eq:model} and
investigates its regularity. Then Section~\ref{sec:4} presents the DG discretization and lists the main results.  Sections~\ref{sec:5} and \ref{sec:6} establish detailed proofs of the error estimates. Finally, Section~\ref{sec:7} conducts several numerical experiments to
verify the theoretical results.

\section{Notation} \label{sec:2}

\setcounter{section}{2} \setcounter{equation}{0}

Firstly, let us make some conventions. For a Lebesgue measurable subset $ \omega $
of $ \mathbb R^l $ ($l=1,2,3$), we use $ H^\gamma(\omega) $
($\gamma\in\mathbb R$) and $ H_0^\gamma(\omega) $ ($ \gamma>0$)
to denote two standard Sobolev spaces \cite{Tartar2007}. For a Lebesgue
measurable subset $ \mathcal O $ of $ \mathbb R^l $ ($ l = 1,2,3,4 $),
the symbol $ \dual{p,q}_{\mathcal O} $ means $ \int_{\mathcal O} pq $ for $pq\in L^1(\mathcal O)$. If $ X $ is a Banach
space, then $ X^* $ denotes its dual space and $ \dual{\cdot, \cdot}_X $ is the
duality pairing between $ X^* $ and $ X $. For $0<\theta<1$ and two Banach spaces $X$ and $Y$, $[X,Y]_{\theta,2}$ stands for the interpolation space constructed via the well-known $K$-method \cite{Tartar2007}, with the norm
\begin{equation}\label{eq:norm-inter}
	\nm{v}_{[X,Y]_{\theta,2}} :=\left(\int_{0}^{\infty}\left(t^{-\theta}K(t,v)\right)^2\frac{\mathrm d t}{t}\right)^{1/2}\quad\forall\,v\in [X,Y]_{\theta,2},
\end{equation}
where the functional $K:(0,\infty)\times (X+Y)\to\mathbb R$ is defined as follows
\[
K(t,v): = \inf_{\substack{v = v_0+v_1\\v_0\in X,\,v_1\in Y}}\left\{
\nm{v_0}_X+t\nm{v_1}_Y
\right\}.
\]
Moreover, if the symbol $ C $ has subscript(s),
then it means a positive constant that depends only on its subscript(s), and its
value may differ at each of its occurrence(s).

Secondly, we introduce some spaces constructed by the eigenvectors of $ -\Delta
$. It is standard that \cite[Theorem 1, \S 6.5.1]{Evansbook} there exists an
orthonormal basis $ \{\phi_n: n \in \mathbb N \} $ of $
L^2(\Omega) $ such that $\phi_n\in H_0^1(\Omega)\cap H^2(\Omega)$ and $ -\Delta \phi_n = \lambda_n \phi_n, $ where $ \{
\lambda_n: n \in \mathbb N \}  $ is a non-decreasing
sequence and $ \lambda_n \to \infty $ as $ n\to \infty $. For any $\gamma\in\mathbb R$, define
\[
\dot H^\gamma(\Omega) := \left\{
\sum_{n=0}^\infty c_n \phi_n:\
\sum_{n=0}^\infty c_n^2 \lambda_n^\gamma < \infty
\right\},
\]
and equip this space with the inner product
\[
\left(
\sum_{n=0}^\infty c_n \phi_n, \sum_{n=0}^\infty d_n \phi_n
\right)_{\dot H^{\gamma(\Omega)}} :=
\sum_{n=0}^\infty \lambda_n^\gamma c_n d_n
\]
for all $ \sum_{n=0}^\infty c_n \phi_n, \sum_{n=0}^\infty d_n \phi_n \in \dot
H^\gamma(\Omega) $, and we use $ \nm{\cdot}_{\dot H^\gamma(\Omega)} $ to denote
the norm induced by this inner product. It is evident that $ \dot
H^\gamma(\Omega) $ is a Hilbert space with an orthonormal basis $
\{\lambda_n^{-\gamma/2} \phi_n: n \in \mathbb N \} $. In addition, $ \dot
H^{-\gamma}(\Omega) $ is the dual space of $ \dot H^\gamma(\Omega) $ in the
sense that
\[
\dual{\sum_{n=0}^\infty
	c_n \phi_n, \sum_{n=0}^\infty d_n \phi_n}_{\dot
	H^\gamma(\Omega)} :=
\sum_{n=0}^\infty c_n d_n
\]
for all $ \sum_{n=0}^\infty c_n \phi_n \in \dot H^{-\gamma}(\Omega) $ and $
\sum_{n=0}^\infty d_n \phi_n \in \dot H^\gamma(\Omega) $.

Thirdly, we introduce some interpolation spaces and vector-valued spaces. In what follows, assume $-\infty<a<b<\infty$. Following \cite{Luo2019}, for any $ m \in\mathbb N $, define
\begin{align*}
	{}^0\!H^m(a,b):={}&\{v\in H^m(a,b): v^{(k)}(b)=0,\,\, 0\leqslant k<m,~k\in\mathbb N\},\\
	{}_0H^m(a,b):={}&\{v\in H^m(a,b): v^{(k)}(a)=0,\,\, 0\leqslant k<m,~k\in\mathbb N\},
\end{align*}
where $ v^{(k)} $ is the $ k $-th weak derivative of $ v $, and endow those two spaces with the following norms
\begin{align*}
	\nm{v}_{{}^0\!H^m(a,b)}:={}&\big\|v^{(m)}\big\|_{L^2(a,b)}\quad\forall\,v\in{}^0\!H^m(a,b),\\
	\nm{v}_{{}_0H^m(a,b)}:={}&\big\|v^{(m)}\big\|_{L^2(a,b)}\quad\forall\,v\in{}_0H^m(a,b).
\end{align*}
Then for $\gamma>0$, we define two interpolation spaces
\[
\begin{split}
	{}^0\!H^{\gamma}(a,b):={}&[L^2(a,b),{}^0\!H^m(a,b)]_{\theta,2},\\
	{}_0H^{\gamma}(a,b):={}&[L^2(a,b),{}_0H^m(a,b)]_{\theta,2},
\end{split}
\]
with corresponding interpolation norms defined by \eqref{eq:norm-inter}, where $0<\theta<1$ and $m\in\mathbb N$ such that $\gamma = m\theta$.
By \cite[Chapter 1]{Lions1972} and \cite{Luo2019}, if $0<\gamma<1/2$, then 	${}_0H^{\gamma}(a,b),\,{}^0\!H^{\gamma}(a,b)$ and $H^{\gamma}(a,b)$ are equivalent and they share the same alternative norm
\begin{equation}\label{eq:norm-1/2}
	\snm{v}_{H^\gamma(a,b)} := \left(
	\int_{\mathbb R} \snm{\xi}^{2\gamma}
	\snm{\mathcal F(Ev)(\xi)}^2 \,\mathrm{d}\xi
	\right)^{1/2} \quad \forall \,v \in H^\gamma(a,b),
\end{equation}
where $ \mathcal F: L^2(\mathbb R) \to L^2(\mathbb R) $ is the Fourier transform and $ Ev$ means extending $v$ to $\mathbb R\backslash(a,b)$ by zero; if $1/2<\gamma<1$, then
\[
\begin{split}
	{}_0H^\gamma(a,b)=\{v\in H^{\gamma}(a,b):
	{}&v^{(m)}(a )= 0\},\\
	{}^0\!H^\gamma(a,b)=\{v\in H^{\gamma}(a,b):
	{}&v^{(m)}(b) = 0\};
\end{split}
\]
if $\gamma = m+s$ with $m\in\mathbb N$ and $s\in(0,1)\backslash\{1/2\}$, then
\[
\begin{split}
	{}_0H^\gamma(a,b)=\Big\{v\in {}_0H^m(a,b):v^{(m)}\in {}_0H^{\gamma-m}(a,b)\Big\},\\
	{}^0\!H^\gamma(a,b)=\Big\{v\in {}^0\!H^m(a,b):v^{(m)}\in {}^0\!H^{\gamma-m}(a,b)\Big\};
\end{split}
\]
if $\gamma = m+1/2$ with $m\in\mathbb N$, then
\[
\begin{split}
	{}_0H^\gamma(a,b)\!=\!{}&\{v\in {}_0H^m(a,b):v^{(m)}\!\in\! H^{1/2}(a,b),\,(t\!-\!a)^{-1/2}v^{(m)}\in L^2(a,b)\},\\
	{}^0\!H^\gamma(a,b)\!=\!{}&\{v\in {}^0\!H^m(a,b):v^{(m)}\!\in\! H^{1/2}(a,b),\, (b\!-\!a)^{-1/2}v^{(m)}\in L^2(a,b)\}.
\end{split}
\]
Moreover, by \cite[Remark 11.5]{Lions1972}, the interpolation norms $\nm{v}_{{}_0H^{1/2}(a,b)}$ and $\nm{v}_{{}^0\!H^{1/2}(a,b)}$ are also equivalent to $\nm{v}_{H^{1/2}(a,b)}+\|(t-a)^{-1/2}v\|_{L^2(a,b)}$ and $\nm{v}_{H^{1/2}(a,b)}+\|(b-t)^{-1/2}v\|_{L^2(a,b)}$, respectively.
Now let $ X $ be a separable Hilbert space with an inner product $ (\cdot,\cdot)_X $ and an
orthonormal basis $ \{e_i: i \in \mathbb N\} $.

For any $\gamma\in\mathbb R$, define the vector-valued space
\[
H^\gamma(a,b;X) := \left\{
v \in L^2(a,b;X):\
\sum_{i=0}^\infty \nm{(v,e_i)_X}_{H^\gamma(a,b)}^2 < \infty
\right\},
\]
and endow this space with the norm
\[
\nm{v}_{H^\gamma(a,b;X)} :=
\left(
\sum_{i=0}^\infty \nm{(v,e_i)_X}_{H^\gamma(a,b)}^2
\right)^{1/2}
\quad \forall\, v \in H^\gamma(a,b;X).
\]
For $\gamma>0$, the two spaces ${}_0H^\gamma(a,b;X)$ and ${}^0\!H^\gamma(a,b;X)$
can be defined analogously.

Fourthly, we introduce the Riemann\textendash Liouville fractional calculus operators. For $ \gamma>0 $, define
\begin{align*}
	\left(\I_{a+}^\gamma v\right)(t) &:= \frac1{ \Gamma(\gamma) }
	\int_a^t (t-s)^{\gamma-1} v(s) \, \mathrm{d}s, \quad t\in(a,b), \\
	\left(\I_{b-}^\gamma v\right)(t) &:= \frac1{ \Gamma(\gamma) }
	\int_t^b (s-t)^{\gamma-1} v(s) \, \mathrm{d}s, \quad t\in(a,b),
\end{align*}
for all $ v \in L^1(a,b;X) $, where $ \Gamma(\cdot) $ denotes the well-known Gamma function. For $
j-1 \leqslant \gamma < j $ with $ j \in \mathbb N_{+} $, define
\[
\D_{a+}^\gamma := \D^j \I_{a+}^{j-\gamma}, \quad
\D_{b-}^\gamma := (-1)^j \D^j \I_{b-}^{j-\gamma},
\]
where $ \D $ is the first-order differential operator in the distribution sense. Essentials of the fractional calculus are listed in Section~\ref{sec:7}.

\section{Weak solution} \label{sec:3}

\setcounter{section}{3} \setcounter{equation}{0}

\subsection{\bf The case $ u_0 = 0 $} \label{subsec:3.1} 
Define
\begin{align*}
	\mathcal X &:= {}_0H^{\alpha/2}(0,T; L^2(\Omega)) \cap L^2(0,T;\dot H^1(\Omega)), \\
	\mathcal Y &:=
	{}^0\!H^{1-\alpha/2}(0,T; L^2(\Omega)) \cap {}^0\!H^{1-\alpha}(0,T;\dot H^1(\Omega)),
\end{align*}
and endow them with the following two norms
\begin{align*}
	\nm{\cdot}_{\mathcal X} &:= \left(
	\nm{\cdot}_{{}_0H^{\alpha/2}(0,T; L^2(\Omega))}^2+
	\nm{\cdot}_{L^2(0,T;\dot H^1(\Omega))}^2
	\right)^{1/2}, \\
	\nm{\cdot}_{\mathcal Y} &:=
	\left(
	\nm{\cdot}_{{}^0\!H^{1-\alpha/2}(0,T; L^2(\Omega))}^2+
	\nm{\cdot}_{{}^0\!H^{1-\alpha}(0,T;\dot H^1(\Omega))}^2
	\right)^{1/2}.
\end{align*}
Assuming that $ f \in \mathcal Y^* $, we call $ u \in \mathcal X $ a weak
solution to problem \eqref{eq:model} if
\begin{equation}
	\label{eq:weak_sol}
	\dual{\D_{0+}^\alpha u, v}_{{}_0H^{\alpha/2}(0,T;L^2(\Omega))} +
	\dual{\nabla u, \nabla v}_{ \Omega_T} =
	\dual{f, \I_{T-}^{1-\alpha} v}_{\mathcal Y}
	\quad  \forall\,v \in \mathcal X.
\end{equation}
Since ${}_0H^{\alpha/2}(0,T;L^2(\Omega))=
H^{\alpha/2}(0,T;L^2(\Omega))$
in the sense of equivalent norms and applying Lemma~\ref{lem:regu} implies that
\begin{equation}\label{eq:injec}
	\nm{\I_{T-}^{1-\alpha} v}_{\mathcal Y}
	\leqslant C_{\alpha} \nm{v}_{\mathcal X}
	\quad \text{ for all } v \in \mathcal X,
\end{equation}
we readily conclude that the above weak solution is well-defined, according to Lemma~\ref{lem:coer}
and the well-known Lax--Milgram theorem.

\begin{theorem}
	\label{thm:basic_weak_sol}
	If $ f \in \mathcal Y^* $, then problem \eqref{eq:model} admits a unique weak
	solution $ u \in \mathcal X$ satisfying $		\nm{u}_{\mathcal X} \leqslant
	C_{\alpha} \nm{f}_{\mathcal Y^*}$.
\end{theorem}

\begin{remark}\label{rem:compare-weak-f}
	Our weak formulation \eqref{eq:weak_sol} is motivated from that proposed in \cite{Li2018A,li_space-time_2009} for \eqref{eq:model-equi}. It looks quite different from the standard one for classical heat equation (cf. \cite[Eq.(2.8)]{schotzau_time_2000} or \cite[Chapter 3]{Lions1972}). A more natural formulation shall be stated as follows : for given $f\in \mathcal Y^*$, find $	u\in \mathcal X$ such that
	\begin{equation}\label{eq:weak-stand}
		-\dual{u,w'}_{\mathcal X}+
		\dual{\nabla u, \D_{T-}^{1-\alpha} \nabla w}_{ \Omega_T}= \dual{f,w}_{\mathcal Y}\quad\forall\,w\in\mathcal Y.
	\end{equation}
	We claim that \eqref{eq:weak-stand} boils down to our nonstandard one \eqref{eq:weak_sol}. Rigorous explanations are provided as follows.
	
	First of all, the operator $\I_{T-}^{1-\alpha}:\mathcal X\to\mathcal Y$ is one-to-one, which means $\mathcal Y = \I_{T-}^{1-\alpha}\mathcal X$. Indeed, from \eqref{eq:injec} follows the injection, and by \cite[Lemmas 3.5 and 3.6]{Luo2019}, $w = \I_{T-}^{1-\alpha}\D_{T-}^{1-\alpha}w$ with $\D_{T-}^{1-\alpha}w\in\mathcal X$, for all $w\in\mathcal Y$.
	
	We then prove that
	\[
	\dual{\D_{0+}^{\alpha}v,\D_{T-}^{1-\alpha}w}_{{}_0H^{\alpha/2}(0,T;L^2(\Omega))}
	=-\dual{w',v}_{{}_0H^{\alpha/2}(0,T;L^2(\Omega))} ,
	\]
	for all $v\in\mathcal X$ and $w\in\mathcal Y$. It is sufficient to consider the scalar case
	\begin{equation}\label{eq:v'w}
		\dual{\D_{0+}^{\alpha/2}v,\D_{T-}^{1-\alpha/2}w}_{(0,T)}
		=-\dual{w',v}_{{}_0H^{\alpha/2}(0,T)} ,
	\end{equation}
	for all $v\in {}_0H^{\alpha/2}(0,T)$ and $w\in {}^0\!H^{1-\alpha/2}(0,T)$. Define a linear functional $s_w:{}_0H^{\alpha/2}(0,T)\to\mathbb R$ by that
	\[
	s_w(v) :=\dual{\D_{0+}^{\alpha/2}v,\D_{T-}^{1-\alpha/2}w}_{(0,T)}
	\quad\forall \,v\in
	{}_0H^{\alpha/2}(0,T),
	\]
	where $w\in{}^0\!H^{1-\alpha/2}(0,T)$ is fixed.
	For any $v\in C_0^\infty(0,T)$, using integration by part, we find $s_w(v) = \dual{v',w}_{(0,T)}$ and thus $s_w = -w'$. Since $C_0^\infty(0,T)$ is dense in ${}_0H^{\alpha/2}(0,T)$, \eqref{eq:v'w} follows immediately.
	
	Hence, by choosing ``indirect" test function $w = \I_{T-}^{1-\alpha}v$ with $v\in \mathcal X$, \eqref{eq:weak-stand} agrees with \eqref{eq:weak_sol}. For investigating the weak solution, it is more convenient for us to treat the latter. This coincides with the same spirit as the weak form proposed for the fractional wave equation in \cite{Luo2019}. However, for numerical discretization, the computational cost of \eqref{eq:weak_sol} is larger than that of \eqref{eq:weak-stand}, since the former has one more fractional integral operator than the latter.
\hfill\ensuremath{\blacksquare}
\end{remark}

Then let us analyze the regularity of the weak solution to problem
\eqref{eq:model}. We first consider the following problem: seek $ y \in
{}_0H^{\alpha/2}(0,T) $ such that
\begin{equation}
	\label{eq:ode}
	\dual{ \D_{0+}^\alpha y, z }_{ {}_0H^{\alpha/2}(0,T) } +
	\lambda \dual{y,z}_{(0,T)} = \dual{ g, \I_{T-}^{1-\alpha} z }_{ {}^0\!H^{1-\alpha/2}(0,T) }
\end{equation}
for all $ z \in {}_0H^{\alpha/2}(0,T) $,
where $ g \in ({}^0\!H^{1-\alpha/2}(0,T))^*$ and $
\lambda >0$ is a constant. Again, by  Lemmas~\ref{lem:coer} and \ref{lem:regu} and the
Lax--Milgram theorem, we conclude that problem \eqref{eq:ode} admits a unique solution
$ y \in {}_0H^{\alpha/2}(0,T) $ and there holds the estimate
\[
\nm{y}_{{}_0H^{\alpha/2}(0,T)} + \lambda^{1/2} \nm{y}_{L^2(0,T)} \leqslant
C_{\alpha} \nm{g}_{({}^0\!H^{1-\alpha/2}(0,T))^*}.
\]

\begin{lemma}
	\label{lem:regu_ode_improved}
	If $ g \in L^2(0,T) $, then the solution $ y $ to problem \eqref{eq:ode} satisfies $		y' + \lambda \D_{0+}^{1-\alpha} y ={} g$ and
	\begin{align}
		\label{eq:regu_ode_est1}
		\nm{y}_{{}_0H^1(0,T)} +
		\lambda \nm{y}_{{}_0H^{1-\alpha}(0,T)}
		\leqslant {}&C_{\alpha} \nm{g}_{L^2(0,T)}.
	\end{align}
	In addition, if $1/2\leqslant \alpha<1$, then for all $0<\epsilon\leqslant 2$,
	\begin{equation}
		\label{eq:regu_ode_C0}
		\lambda^{1/(2\alpha)-\sigma\epsilon/2}\nm{y}_{C[0,T]}\leqslant
		\frac{		C_{\alpha,T} }{\epsilon^{\sigma/2}}\nm{g}_{L^2(0,T)},
	\end{equation}
	where $\sigma = 0$ if $1/2<\alpha<1$ and $\sigma=1$ if $\alpha=1/2$.
\end{lemma}

\begin{proof}
Observe that the equality
\begin{equation}
	\label{eq:ode-y}
	y = \I_{0+}^\alpha(\I_{0+}^{1-\alpha} g - \lambda y) =
	\I_{0+} g - \lambda \I_{0+}^\alpha y
\end{equation}
is contained in the proofs of \cite[Lemmas 3.1 and 3.2]{Li2018A},
and using \eqref{eq:ode-y} recursively yields the relation
\[
\begin{split}
	y
	={}&(-\lambda)^{n}\I_{0+}^{n\alpha}y+ \sum_{i=0}^{n-1}(-\lambda)^{i}\I_{0+}^{1+i\alpha}g
\end{split}
\qquad \forall\,
n \in \mathbb N.
\]
Then by Lemma~\ref{lem:regu}, $ y \in {}_0H^1(0,T) $,
and it follows from \eqref{eq:ode-y} that
\begin{equation}
	\label{eq:regu_ode_strong}
	y' + \lambda \D_{0+}^{1-\alpha} y ={}g.
\end{equation}

Next, let us prove \eqref{eq:regu_ode_est1}. Multiplying both sides of
\eqref{eq:regu_ode_strong} by $ y' $ and integrating over $ (0,T) $ gives
\[
\nm{y'}_{L^2(0,T)}^2 + \lambda \dual{\D_{0+}^{1-\alpha} y, y'}_{(0,T)}
= \dual{g, y'}_{(0,T)},
\]
so that, by Lemmas~\ref{lem:coer} and \ref{lem:regu}, the Cauchy--Schwartz inequality and Young's
inequality with $ \epsilon $, we obtain the estimate
\begin{equation}
	\label{eq:regu-1}
	\nm{y'}_{L^2(0,T)}^2 + \lambda
	\nm{y}_{{}_0H^{1-\alpha/2}(0,T)}^2
	\leqslant C_{\alpha} \nm{g}_{L^2(0,T)}^2.
\end{equation}
Additionally, since $ y \in {}_0H^1(0,T) $, a straightforward
computing gives
\begin{equation}\label{eq:IDy}
	y = \I_{0+}^{1-\alpha}\D_{0+}^{1-\alpha}y,
\end{equation}
which together with \eqref{eq:regu_ode_strong} and Lemma~\ref{lem:regu} implies
\begin{equation}
	\label{eq:regu-2}
	\lambda \nm{y}_{{}_0H^{1-\alpha}(0,T)} \leqslant
	C_{\alpha} \left( \nm{y'}_{L^2(0,T)} + \nm{g}_{L^2(0,T)} \right).
\end{equation}
Therefore, combining
\eqref{eq:regu-1} and \eqref{eq:regu-2} proves \eqref{eq:regu_ode_est1}.

It remains to prove \eqref{eq:regu_ode_C0}. Thanks to \eqref{eq:regu_ode_est1} and Lemma~\ref{lem:Nirenberg}, for $1/2<\alpha<1$, we have
\begin{align*}
	\lambda^{1/(2\alpha)} \nm{y}_{C[0,T]}
	\leqslant{} 
	C_{\alpha,T}
	\nm{y}_{{}_0H^{1}(0,T)}^{(\alpha-1/2)/\alpha}
	\big(
	\lambda \nm{y}_{{}_0H^{1-\alpha}(0,T)}
	\big)^{1/(2\alpha)} 
	\leqslant {} C_{\alpha,T}
	\nm{g}_{L^{2}(0,T)}.
\end{align*}
As for $\alpha=1/2$, applying Lemma~\ref{lem:Nirenberg} once again yields that
\begin{align*}
	\lambda^{1-\epsilon/2} \nm{y}_{C[0,T]} \leqslant {}
	\frac{C_{\alpha,T}}{\sqrt{\epsilon}}
	\big(
	\lambda \nm{y}_{{}_0H^{1/2}(0,T)}
	\big)^{1-\epsilon/2}
	\nm{y}_{{}_0H^{1}(0,T)}^{ \epsilon/2}
	\leqslant {}
	\frac{C_{\alpha,T}}{\sqrt{\epsilon}}
	\nm{g}_{L^{2}(0,T)}
\end{align*}
for all $0<\epsilon\leqslant2$. This completes the proof of this lemma.
\end{proof}

\begin{lemma}
	\label{lem:regu_ode_higher}
	If $ g \in {}_0H^{\gamma}(0,T) $ with $0<\gamma\leqslant1/2$,
	then the solution $ y $ to problem \eqref{eq:ode} satisfies
	that
	\begin{equation}
		\label{eq:regu_ode_est-high}
		\nm{y}_{{}_0H^{1+\gamma}(0,T)} +
		\lambda \nm{y}_{{}_0H^{1+\gamma-\alpha}(0,T)}
		\leqslant {}C_{\alpha,\gamma} \nm{g}_{{}_0H^{\gamma}(0,T)}.
	\end{equation}
\end{lemma}
\begin{proof}
Let us focus on the auxiliary problem: seek $ w\in
{}_0H^{\alpha/2}(0,T) $ such that
\begin{equation}
	\label{eq:ode-w}
	\dual{ \D_{0+}^\alpha w, z }_{ {}_0H^{\alpha/2}(0,T) } +
	\lambda \dual{w,z}_{(0,T)} = \dual{ \D_{0+}^{\gamma}g, \I_{T-}^{1-\alpha} z }_{(0,T) }
\end{equation}
for all $ z \in {}_0H^{\alpha/2}(0,T) $. By Lemmas~\ref{lem:regu_ode_improved} and \ref{lem:coer1}, $ w \in {}_0H^{1}(0,T) $ exists uniquely and satisfies
\begin{equation}\label{eq:ode-strang-w}
	w'+ \lambda \D_{0+}^{1-\alpha} w = \D_{0+}^{\gamma} g
\end{equation}
In addition, by \eqref{eq:regu_ode_est1} and Lemma~\ref{lem:coer1}, it holds that
\begin{equation*} 	
	\nm{w}_{{}_0H^{1}(0,T)} +
	\lambda \nm{w}_{{}_0H^{1-\alpha}(0,T)}
	\leqslant {}C_{\alpha}
	\nm{\D_{0+}^{\gamma}g}_{L^{2}(0,T)}
	\leqslant {}C_{\alpha,\gamma} \nm{g}_{{}_0H^{\gamma}(0,T)}.
\end{equation*}

Set $y:=\I_{0+}^\gamma w$, then using Lemma~\ref{lem:regu} gives \eqref{eq:regu_ode_est-high}. It remains to prove that
$y=\I_{0+}^\gamma w$ is the solution to problem \eqref{eq:ode}. To this end, applying $\I_{0+}^\gamma$ to both sides of \eqref{eq:ode-strang-w}, we obtain
\[
(\I_{0+}^\gamma w)'+ \lambda \D_{0+}^{1-\alpha} \I_{0+}^\gamma w= g.
\]
This proves the desired result and completes the proof.
\end{proof}

If $ f \in L^2(0,T;\dot H^{-1}(\Omega)) $, then similar to \cite[Lemma 4.4]{Luo2019} and \cite[Theorem 3.1]{Li2018A}, we can prove that the weak solution $ u $ to problem \eqref{eq:model} is
\[
u(t) = \sum_{n=0}^\infty y_n(t) \phi_n, \quad 0 < t \leqslant  T,\]
where $ y_n \in {}_0H^{\alpha/2}(0,T) $ satisfies that
\[
\dual{\D_{0+}^\alpha y_n, z}_{{}_0H^{\alpha/2}(0,T)} +
\lambda_n \dual{y_n, z}_{(0,T)} =
\dual{\dual{f,\phi_n}_{\dot H^1(\Omega)}, \I_{T-}^{1-\alpha} z}_{(0,T)}
\]
for all $ z \in {}_0H^{\alpha/2}(0,T) $. Therefore, the desired regularity result follows directly from Lemmas~\ref{lem:regu_ode_improved} and \ref{lem:regu_ode_higher}.
\begin{theorem}
	\label{thm:regu_pde}
	Assume that $ f \in {}_0H^{\gamma}(0,T;\dot H^{-\beta}(\Omega)) $ with $0\leqslant \gamma\leqslant 1/2$ and $ 0 \leqslant \beta
	\leqslant 1 $. Then the weak solution $ u $ to problem \eqref{eq:model} satisfies
	$u' - \D_{0+}^{1-\alpha} \Delta u ={} f$ in $L^2(0,T; \dot{H}^{-\beta}(\Omega))$
	and
	\begin{align*}
		{}\nm{u}_{{}_0H^{1+\gamma}(0,T;\dot H^{-\beta}(\Omega))} +
		\nm{u}_{{}_0H^{1+\gamma-\alpha}(0,T;\dot H^{2-\beta}(\Omega))}
		\leqslant{} &
		C_{\alpha,\gamma} \nm{f}_{{}_0H^{\gamma}(0,T;\dot H^{-\beta}(\Omega))}.
	\end{align*}
	In addition, if $\gamma=0$ and $1/2\leqslant \alpha<1$, then for all $0<\epsilon\leqslant 2$,
	\begin{equation*}
		\nm{y}_{C([0,T];\dot{H}^{1/\alpha-\sigma\epsilon-\beta}(\Omega))}
		\leqslant
		\frac{		C_{\alpha,T} }{\epsilon^{\sigma/2}}\nm{f}_{L^2(0,T;\dot H^{-\beta}(\Omega))},
	\end{equation*}
	where $\sigma = 0$ if $1/2<\alpha<1$ and $\sigma=1$ if $\alpha=1/2$.
\end{theorem}

For the dual problem of \eqref{eq:model}, we have the following theorem.
\begin{theorem}
	\label{thm:dual}
	Assume that $ q \in L^2(0,T;\dot H^{-\beta}(\Omega)) $ with $ 0 \leqslant \beta
	\leqslant 1 $. Then there exists a unique
	\[
	w\in \mathcal G:={{}^0\!}H^{1}(0,T;\dot H^{-\beta}(\Omega)) \cap {{}^0\!}H^{1-\alpha}(0,T;\dot H^{2-\beta}(\Omega))
	\]
	such that $		-w' - \D_{T-}^{1-\alpha} \Delta w ={} q$ and
	\begin{align*}
		{}\nm{w}_{{{}^0\!}H^{1}(0,T;\dot H^{-\beta}(\Omega))} +
		\nm{w}_{{{}^0\!}H^{1-\alpha}(0,T;\dot H^{2-\beta}(\Omega))}
		\leqslant{} &
		C_{\alpha} \nm{q}_{L^{2}(0,T;\dot H^{-\beta}(\Omega))}.
	\end{align*}
\end{theorem}

\subsection{\bf The case $ f=0 $} \label{subsec:3.2} 

For $a,b>0$, recall the Mittag-Leffler function
\[
E_{a,b}(z) := \sum_{k=0}^\infty
\frac{z^k}{\Gamma(ak + b)},
\quad z \in \mathbb C.
\]
Given $ \lambda, t > 0$ and $\gamma\in\mathbb R_+\backslash \mathbb N$, we have the following facts
(cf. \cite{gorenflo_mittag-leffler_2014})： 
\begin{align}
	\snm{E_{a,b}(-t)} \leqslant
	\frac{C_{a,b}}{1+t},
	\label{eq:Eabz-est}
	\\
	\D_{0+}^\gamma E_{a,1}(-\lambda t^a) =
	t^{-\gamma}E_{a,1-\gamma}(-\lambda t^a),
	\label{eq:DgE}
	\\
	\frac{\mathrm{d}}{\mathrm d t}\, E_{a,1}(-\lambda t^a) =
	-\lambda t^{a-1}E_{a,a}(-\lambda t^a).
	\label{eq:DE}
\end{align}

For any $\lambda>0$ and $y_0\in\mathbb R$, by \eqref{eq:DgE} and \eqref{eq:DE}, it is easy to see that
\[
y(t) =  y_0E_{\alpha,1}(-\lambda t^{\alpha}),\quad 0\leqslant t\leqslant T,
\]
solves the equation
\[
y' +\lambda \D_{0+}^{1-\alpha} y = 0,\quad 0<t\leqslant T,
\]
with initial condition $y(0) = y_0$. Therefore, for $f = 0$ and $u_0\in\dot{H}^{-2}(\Omega)$, it is natural to define a weak solution
of problem \eqref{eq:model} by that \cite{sakamoto_initial_2011}： 
\begin{equation}
	\label{eq:weak-solu-homo}
	u(t) : =
	\sum_{n=0}^\infty
	E_{\alpha,1}(-\lambda_n t^{\alpha})
	\dual{u_{0},\phi_n}_{\dot H^2(\Omega)}\phi_n,
	\quad 0 \leqslant t \leqslant T.
\end{equation}
It follows from \eqref{eq:Eabz-est} that
$ u \in C([0,T];\dot H^{-2}(\Omega)) $ is well defined. In addition, we have $u(0) = u_0$ and
\[
\nm{u}_{C([0,T];\dot H^{-2}(\Omega))}
\leqslant C_{\alpha}\nm{u_0}_{\dot H^{-2}(\Omega)}.
\]

Since $u_0\in\dot H^{-2}(\Omega)$, \eqref{eq:weak-solu-homo} shall be understood as the ``very weak solution" by using the transposition method \cite{Lions1972}; see also \cite[Section 4.3]{Luo2019}. In the following, we only consider the case $u_0\in L^2(\Omega)$ and establish the weak formulation \eqref{eq:weak-form-fpde-u0} and regularity estimate \eqref{eq:regu-fpde-u0}. Particularly, the formulation \eqref{eq:weak-form-fpde-u0} tells us in what sense \eqref{eq:weak-solu-homo} is a weak solution to the original subdiffusion model \eqref{eq:model}.
\begin{theorem}
	\label{thm:regu-fpde-u0}
	If $ u_0 \in L^2(\Omega)$,
	then the weak solution defined by \eqref{eq:weak-solu-homo} satisfies
	\begin{equation}
		\label{eq:weak-form-fpde-u0}
		\dual{u', v}_{H^{(1-\alpha)/2}(0,T;\dot H^1(\Omega))} +
		\dual{\D_{0+}^{1-\alpha}\nabla u, \nabla v}_{H^{(1-\alpha)/2}(0,T;L^2(\Omega))} = 0,
	\end{equation}
	for all $v \in H^{(1-\alpha)/2}(0,T;\dot H^1(\Omega))$, and we have the estimate
	\begin{equation}
		\label{eq:regu-fpde-u0}
		\begin{split}
			{}&\nm{u'}_{(H^{(1-\alpha)/2}(0,T;\dot H^{1}(\Omega)))^*}+
			\nm{u}_{C([0,T];L^{2}(\Omega))}
			+\nm{u}_{H^{(1-\alpha)/2}(0,T;\dot H^{1}(\Omega))}
			\\
			{}&\qquad \quad
			+\epsilon_{\alpha,\gamma}
			\nm{u}_{L^{2}(0,T;\dot H^{\gamma}(\Omega))}
			\leqslant {}C_{\alpha,T} \nm{u_0}_{L^2(\Omega)},
		\end{split}
	\end{equation}
	where $\epsilon_{\alpha,\gamma} := \sqrt{2-\gamma}+\sqrt{\snm{2\alpha-1}}$ with $\gamma=\min\{2,1/\alpha\}$ if $\alpha\neq1/2$ and $1\leqslant \gamma<2$ if $\alpha=1/2$.
\end{theorem}

\begin{proof}
By \eqref{eq:DgE} and \eqref{eq:DE}, a routine computation yields
\begin{equation*}
	t^{1-\alpha} \nm{\D_{0+}^{1-\alpha}\!u(t)}_{\dot H^1(\Omega)}^2
	+	t^{1-\alpha}\nm{u'(t)}_{\dot H^{-1}(\Omega)}^2
	\!	\leqslant\!C_\alpha\!\sum_{n=0}^\infty
	\frac{\lambda_nt^{\alpha-1}}
	{(1+\lambda_nt^{\alpha})^2}
	\dual{u_0,\phi_n}_\Omega^2,
\end{equation*}
and
\begin{equation*}
	\nm{\D_{0+}^{(1-\alpha)/2}\!u(t)}_{\dot H^1(\Omega)}^2+
	\nm{\I_{0+}^{(1-\alpha)/2}\!u'(t)}_{\dot H^{-1}(\Omega)}^2
	\!	\leqslant\!C_\alpha \!\sum_{n=0}^\infty
	\frac{\lambda_nt^{\alpha-1}}
	{(1+\lambda_nt^{\alpha})^2}\!
	\dual{u_0,\phi_n}_\Omega^2.
\end{equation*}
Observing Lemma~\ref{lem:coer1} and the estimate
\[
\begin{split}
	{}\int_0^T\frac{\lambda_nt^{\alpha-1}}{(1+\lambda_nt^{\alpha})^2}
	\mathrm d t
	\leqslant {}
	\int_0^{\lambda_n^{-1/\alpha}}
	\lambda_nt^{\alpha-1}
	\mathrm d t +
	\int_{\lambda_n^{-1/\alpha}}^\infty
	\lambda_n^{-1} t^{-1-\alpha}
	\mathrm d t={}\frac{2}{\alpha},
\end{split}
\]
we have that
\[
\begin{split}
	\nm{u}_{H^{(1-\alpha)/2}(0,T;\dot H^{1}(\Omega))}^2
	+\int_{0}^{T}	t^{1-\alpha} \nm{\D_{0+}^{1-\alpha}\nabla u(t)}_{L^2(\Omega)}^2\mathrm d t		{}&\\
	+	\int_{0}^{T}t^{1-\alpha}\nm{u'(t)}_{\dot H^{-1}(\Omega)}^2\mathrm d t
	{}&		\leqslant C_{\alpha} \nm{u_0}_{L^2(\Omega)}^2.
\end{split}
\]
Since using  \cite[Lemma 16.3]{Tartar2007} implies
\[
\int_{0}^{T}	t^{\alpha-1} \nm{v(t)}_{\dot H^1(\Omega)}^2\mathrm d t\leqslant  C_{\alpha,T}
\nm{v}_{H^{(1-\alpha)/2}(0,T;\dot H^1(\Omega))}^2,
\]
for all $	v\in H^{(1-\alpha)/2}(0,T;\dot H^1(\Omega))$, it follows
\begin{equation*} 		
	\nm{u'}_{(H^{(1-\alpha)/2}(0,T;\dot H^{1}(\Omega)))^*}
	+\nm{\D_{0+}^{1-\alpha}\nabla u}_{(H^{(1-\alpha)/2}(0,T;L^{2}(\Omega)))^*}
	\leqslant {}C_{\alpha,T} \nm{u_0}_{L^2(\Omega)},
\end{equation*}
and a direct calculation gives \eqref{eq:weak-form-fpde-u0}. Meanwhile, a similar manipulation accomplishes the estimate \eqref{eq:regu-fpde-u0} and finishes the proof of this theorem.
\end{proof}

\section{Main results} \label{sec:4}

\setcounter{section}{4} \setcounter{equation}{0}

Let us introduce the piecewise constant DG method proposed in
\cite{McLean2009Convergence}. Given $ J \in \mathbb N_{>0} $, let
$0=t_0<t_1<\cdots<t_J=T$ be a partition of $[0,T]$ with $\tau: =\max_{1\leqslant  j\leqslant J}(t_j-t_{j-1})$,
and set $ I_j := (t_{j-1},t_j) $ for $ 1
\leqslant j \leqslant J $. Let $ \mathcal K_h $ be a conventional conforming and
quasi-uniform triangulation of $ \Omega $ consisting of $ d $-simplexes, and we use $
h $ to denote the maximum diameter of the elements in $ \mathcal K_h $. Define
\begin{align*}
	S_h&:= \left\{
	v_h \in H_0^1(\Omega):\
	v_h|_K \in P_1(K) \quad\forall \,K \in \mathcal K_h
	\right\},\\
	\mathcal X_{\tau,h} &:= \left\{
	V \in L^2(0,T; S_h):\
	V|_{I_j} \in P_0(I_j; S_h)\quad\forall\,1\leqslant j\leqslant J
	\right\},
\end{align*}
where $ P_1(K) $ is the set of all linear polynomials defined on $ K $, and $
P_0(I_j;  S_h) $ is the set of all $  S_h$-valued
constant functions on $ I_j $. For each $ V \in \mathcal X_{\tau,h} $, we will
use the following notations:
\begin{align*}
	V_j^{+}  :={}\lim_{t \to {{t_j}+}} V(t)
	\quad\text{for }&0 \leqslant j < J,\text{ and } V_J^+ := 0; \\
	V_j^{-}  := {}\lim_{t \to {{t_j}-}} V({t})
	\quad\text{for }& 1 \leqslant j \leqslant J, \text{ and } V_0^- := 0; \\
	\jmp{V_j}  :={} V_j^{+} - V_j^{-}
	\quad\text{for }& 0 \leqslant j \leqslant  J.
\end{align*}

Assuming that $ u_0 \in  S_h^* $ and $ f \in \mathcal X_{\tau,h}^* $, the piecewise constant DG method defines a numerical solution $ U \in \mathcal X_{\tau,h} $ to problem \eqref{eq:model} by that
\begin{equation}
	\label{eq:numer_sol}
	\mathcal A(U,V)=
	\dual{f,V}_{\mathcal X_{\tau,h}} +
	\dual{u_0,V_0^{+}}_{ S_h}
	\quad \forall\,V \in \mathcal X_{\tau,h},
\end{equation}
where
\[
\mathcal A(W,V) := \sum_{j=0}^{J-1} \dual{\jmp{W_j}, V_j^{+}}_\Omega +
\dual{ \D_{0+}^{1-\alpha} \nabla W, \nabla V }_{\Omega_T}
\]
for all $W,V\in\mathcal X_{\tau,h}$. It is not hard to conclude from \cite[Theorem 12.1]{Thomee2006} and Lemma~\ref{lem:coer} that
if $ V \in \mathcal  X_{\tau,h} $, then we have
\begin{equation}\label{eq:sj}
	\begin{aligned}
		\mathcal A(V,V\chi_{(0,t_j)})\geqslant {}&
		\frac12 \left(
		\nm{V_j^{-}}_{L^2(\Omega)}^2 + \nm{V_0^{+}}_{L^2(\Omega)}^2\right)\\
		{}&\quad+\sin
		\frac{\alpha\pi}{2}
		\snm{V}_{H^{(1-\alpha)/2}(0,t_j; \dot H^1(\Omega))}^2,
	\end{aligned}
\end{equation}
for all $ 1 \leqslant j \leqslant J $. Above and in what follows, $ \chi_{(a,b)} $ denotes the indicator function of the interval $ (a,b) $.

For convenience, in what follows we assume that $ u $ is the weak solution to problem
\eqref{eq:model} and $ U $ is its numerical approximation defined by
\eqref{eq:numer_sol}. The notation $ a \lesssim b $ means that there exists a
generic positive constant $ C $, independent of $ h $, $ \tau $ and $ u $, such
that $ a \leqslant C b $. Moreover, $ a \sim b $ means $ a \lesssim b \lesssim a$.

The well-posedness of the solution $U$ to \eqref{eq:numer_sol} is firstly established in the following theorem.

\begin{theorem}
	\label{thm:stab}
	If $ u_0 \in L^2(\Omega) $ and $ f \in (H^{(1-\alpha)/2}(0,T;\dot
	H^{1}(\Omega)) )^*$, then problem \eqref{eq:numer_sol} admits a unique solution $U$ such that
	\begin{equation}
		\label{eq:stab}
		\begin{aligned}
			{}			&
			\nm{U}_{L^\infty(0,T;L^2(\Omega))}+
			\snm{U}_{H^{(1-\alpha)/2}(0,T;\dot H^1(\Omega))} \\
			\lesssim{}&	\nm{u_0}_{L^2(\Omega)}+
			\nm{f}_{(H^{(1-\alpha)/2}(0,T;\dot
				H^{1}(\Omega)) )^*}.
		\end{aligned}
	\end{equation}
\end{theorem}

\begin{remark}
	We also refer the reader to \cite[Theorem 1]{McLean2009Convergence} for
	another stability estimate, which is derived in the case that $ u_0 \in
	L^2(\Omega) $ and $ f \in L^1(0,T;L^2(\Omega)) $.\hfill\ensuremath{\blacksquare}
\end{remark}

Below, let us present our main error estimates.

\begin{theorem}
	\label{thm:conv-f-L2}
	If $ u_0 = 0 $ and $ f \in L^2(\Omega_T) $, then
	\begin{align}
		\label{eq:conv-f-L2-L2}
		\nm{u-U}_{L^2(\Omega_T)}    \lesssim{} &
		\big( h^2 + \tau \big)
		\nm{f}_{L^2(\Omega_T)},\\
		\label{eq:conv-f-L2-Hs}
		\snm{u-U}_{H^{(1-\alpha)/2}(0,T;\dot H^1(\Omega))}    \lesssim{} &
		\big( h + \tau^{1/2} \big)
		\nm{f}_{L^2(\Omega_T)}.
	\end{align}
\end{theorem}

\begin{theorem}
	\label{thm:conv-u0-L2}
	Assume that $ f = 0 $. If $ u_0 \in L^2(\Omega) $, then
	\begin{equation}
		\label{eq:conv-u0-L2-L2}
		\nm{u-U}_{L^2(\Omega_T)}
		\lesssim \big( h + \tau^{1/2} \big) \nm{u_0}_{L^2(\Omega)}.
	\end{equation}
\end{theorem}

\begin{remark}\label{rem:opt-rate}
	In view of Lemma~\ref{lem:interp}, Theorems~\ref{thm:regu_pde} and \ref{thm:regu-fpde-u0}, we conclude that all the convergence rates in Theorem~\ref{thm:conv-f-L2} are optimal with respect to the solution regularity while Theorem ~\ref{thm:conv-u0-L2} only gives optimal estimate in temporal discretization. Indeed, ignoring the logarithm factor, the spatial accuracy in \eqref{eq:conv-u0-L2-L2} should be $\min\{2,1/\alpha\}$; see numerical results in Table~\ref{tab:subdf-ex3-space}.\hfill\ensuremath{\blacksquare}
\end{remark}

\begin{remark}
	Although the half order $O(\tau^{1/2})$ in \eqref{eq:conv-u0-L2-L2} is optimal with respect to the solution regularity (cf. Theorem~\ref{thm:regu-fpde-u0}) and this has been verified by the numerical test with uniform temporal grid (see Table~\ref{tab:subdf-ex3-temp}), it is still possible to recover the first order accuracy $O(\tau)$ by using graded temporal meshes; see \cite[Theorem 3.1]{Li2020Numer} for rigorous proof of the L1 scheme with smoother data $u_0\in \dot H^\sigma(\Omega),\, 0<\sigma\leqslant 2$. In the last experiment of Section~\ref{sec:7}, we investigate the performance of our DG scheme under nonuniform grids and obtain the rate $O(\tau)$ for $u_0\in L^2(\Omega)$ with suitable graded meshes (cf. Table~\ref{tab:ex4-u0}).\hfill\ensuremath{\blacksquare}
\end{remark}

Moreover, if the temporal grid is equi-distributed, then quasi-optimal (including logarithm factors) error bounds under the $L^\infty(0,T;L^2(\Omega))$-norm  are derived.

\begin{theorem}
	\label{thm:conv-Linf}
	Assume $u_0=0$ and the temporal grid is uniform.
	If $ f \in L^2(\Omega_T) $, then
	\begin{align}
		\label{eq:conv-LinfL2-f-L2}
		\nm{u - U}_{L^\infty(0,T;L^{2}(\Omega))}
		\lesssim{}&\snm{\ln \tau}
		\Big(\tau^{1/2} + \epsilon_h
		h^{\min\{2,1/\alpha\}}\Big)
		\nm{f}_{L^2(\Omega_T)},
	\end{align}
	where $\epsilon_h = 1$ if $\alpha\neq1/2$ and $\epsilon_h = \sqrt{\snm{\ln h}}$ if $\alpha=1/2$. Moreover, if $ f \in {}_0H^{1/2}(0,T;L^2(\Omega)) $, then
	\begin{align}
		\label{eq:conv-LinfL2-f-Hs}
		\nm{u - U}_{L^\infty(0,T;L^2(\Omega))}
		\lesssim{}&\snm{\ln \tau}
		\left(\snm{\ln \tau}\tau+  h^2 \right)
		\nm{f}_{{}_0H^{1/2}(0,T;L^2(\Omega))}.
	\end{align}
\end{theorem}


\section{Proofs of Theorems~\ref{thm:stab}\textendash \ref{thm:conv-u0-L2}} \label{sec:5}

\setcounter{section}{5} \setcounter{equation}{0}

\subsection{\bf Preliminaries} \label{subsec:5.1} 

Given a Banach space $ X $, we introduce two interpolation operators as
follows \cite[Chapter~12]{Thomee2006}: given $ v \in C((0,T];X) $ and $ w \in C([0,T);X) $, define $ P_\tau v $ and $ Q_\tau w$ respectively by that
\[
\begin{aligned}
	\left(  P_\tau v \right)|_{I_j} = {}&v(t_j)
	&& \forall\,1 \leqslant j \leqslant J,\\
	\left(  Q_\tau w \right)|_{I_j} = {}&w(t_{j-1})
	&& \forall\,1 \leqslant j \leqslant J.
\end{aligned}
\]
Let $ P_h: L^2(\Omega) \to  S_h $ be an $ L^2(\Omega) $-orthogonal projection operator and $R_h:\dot{H}^1(\Omega)\to S_h$ be the Ritz projection operator.
Then by the theory of interpolation spaces \cite{Tartar2007} and the
standard approximation estimates,
we readily obtain that if $ v \in \dot H^r(\Omega) $ with $
1 \leqslant r \leqslant 2 $, then
\[
\begin{split}
	\nm{(I-R_h)v}_{L^2(\Omega)} +h\nm{(I-R_h)v}_{\dot{H}^1(\Omega)}
	\lesssim{}&
	h^r \nm{v}_{\dot H^r(\Omega)},\\
	\nm{(I-P_h)v}_{L^2(\Omega)} +h\nm{(I-P_h)v}_{\dot{H}^1(\Omega)}
	\lesssim{}&
	h^r \nm{v}_{\dot H^r(\Omega)}.
\end{split}
\]
Because the above two estimates are well known, we will use them implicitly for clarity.

Below, let us establish some nonstandard error estimates of $P_\tau$ and $Q_\tau$.

\begin{lemma}
	\label{lem:PQ-Hs}
	If $ 0 \leqslant \beta < 1/2 < \gamma \leqslant 1 $ and $ v \in H^\gamma(0,T) $, then
	\begin{align}
		\nm{\left( I-P_\tau \right)v}_{H^\beta(0,T)} \!+\!\nm{\left( I-Q_\tau \right)v}_{H^\beta(0,T)}
		\leqslant{} &C_{\beta,\gamma,T} \tau^{\gamma-\beta} \nm{v}_{H^\gamma(0,T)}.
		\label{eq:PQ-Hs}
	\end{align}
\end{lemma}

\begin{proof}
We first consider the estimate for $P_\tau$ and
set $ g := v-P_\tau v $. In view of the proof of \cite[Lemma 4.3]{li_space-time_2019}, we have
\begin{equation}\label{eq:bd-g}
	\nm{g}_{H^\beta(0,T)}^2 \leqslant C_{\beta,T}
	\big( \mathbb I_1 + \mathbb I_2 \big),
\end{equation}
where
\[
\left\{
\begin{aligned}
	\mathbb I_1 &:=  \sum_{j=1}^J
	\int_{I_j} \int_{I_j} \frac{\snm{v(s)-v(t)}^2}{\snm{s-t}^{1+2\beta}}
	\,\mathrm d s\mathrm d t,\\
	\mathbb I_2 &:= 	\sum_{j=1}^J \int_{I_j} g^2(t) \left(
	(t_j-t)^{-2\beta} + (t-t_{j-1})^{-2\beta}
	\right) \mathrm d t .
\end{aligned}
\right.
\]

If $\gamma=1$, then $v\in H^1(0,T)$ and
\begin{align*}
	\mathbb I_1={}&\sum_{j=1}^J
	\int_{I_j} \int_{I_j} \frac{\snm{v(s)-v(t)}^2}{\snm{s-t}^{1+2\beta}}
	\,\mathrm d s\mathrm d t
	={}\sum_{j=1}^J
	\int_{I_j} \int_{I_j} \frac{\snm{\int_{s}^{t}v'(r)
			\mathrm d r}^2}{\snm{s-t}^{1+2\beta}}
	\,\mathrm d s\mathrm d t\\
	\leqslant {}&\sum_{j=1}^J\left(\int_{I_j}\snm{v'(r)}^2\mathrm dr\right)
	\left(\int_{I_j} \int_{I_j} \snm{s-t}^{-2\beta}
	\mathrm d s\mathrm d t\right)\\
	\leqslant{}&C_\beta \sum_{j=1}^J\tau_j^{2-2\beta}
	\int_{I_j}\snm{v'(r)}^2\mathrm dr\leqslant C_\beta \tau^{2-2\beta}\nm{v}_{H^1(0,T)}^2.
\end{align*}
In addition, the term $\mathbb I_2$
can be estimated similarly
\begin{align*}
	\mathbb I_2 ={}& 	\sum_{j=1}^J \int_{I_j} g^2(t) \left(
	(t_j-t)^{-2\beta} + (t-t_{j-1})^{-2\beta}
	\right)\mathrm d t\\
	={}&\sum_{j=1}^J \int_{I_j}  \snm{\int_{t_{j}}^t v'(r)  \, \mathrm{d}r}^2
	\left((t_j-t)^{-2\beta} + (t-t_{j-1})^{-2\beta}
	\right)  \mathrm d t\\
	\leqslant {}&C_\beta \sum_{j=1}^J\tau_j^{2-2\beta}
	\int_{I_j}\snm{v'(r)}^2\mathrm dr\leqslant C_\beta \tau^{2-2\beta}\nm{v}_{H^1(0,T)}^2.
\end{align*}
Plugging the above two estimates into \eqref{eq:bd-g} gives
\[
\nm{\left( I-P_\tau \right)v}_{H^\beta(0,T)}
\leqslant  C_{\beta,T} \tau^{1-\beta}\nm{v}_{H^1(0,T)}.
\]

Next, we consider $1/2<\gamma<1$. Observing \cite[Lemma 36.1]{Tartar2007}, we find that
\[
\mathbb I_1 \leqslant
\sum_{j=1}^J \tau_j^{2(\gamma-\beta)} \int_{I_j} \int_{I_j}
\frac{\snm{v(s)-v(t)}^2}{\snm{s-t}^{1+2\gamma}} \,\mathrm d s\mathrm d t
\leqslant C_{\gamma,T}
\tau^{2(\gamma-\beta)} \nm{v}_{H^\gamma(0,T)}^2.
\]
We aim to establish
\begin{equation}\label{eq:I2}
	\mathbb I_2 \leqslant
	C_{\gamma,T} \tau^{2(\gamma-\beta)} \nm{v}_{H^\gamma(0,T)}^2.
\end{equation}
Since $ 0 \leqslant \beta < 1/2 < \gamma <1 $, by \cite[Theorem 5.20]{Kufner2017}, we have
\[
\begin{aligned}
	{}&	\int_0^1 |w(t)-w(1)|^2
	\left( t^{-2\beta} + (1-t)^{-2\beta} \right)\mathrm d t\\	
	\leqslant 		{}&\int_0^1 |w(t)-w(1)|^2
	\left( t^{-2\gamma} +(1-t)^{-2\gamma} \right)\mathrm d t\\
	\leqslant{}&	C_{\gamma}
	\int_0^1 \int_0^1
	\frac{\snm{w(s)-w(t)}^2}{\snm{s-t}^{1+2\gamma}} \,\mathrm d s\mathrm d t,
\end{aligned}
\]
for all $ w \in H^\gamma(0,1) $.
Therefore, a standard scaling argument implies
\begin{align*}
	\mathbb I_2 & \leqslant
	C_{\gamma} \sum_{j=1}^J
	\tau_j^{2(\gamma-\beta)} \int_{I_j} \int_{I_j}
	\frac{\snm{v(s)-v(t)}^2}{\snm{s-t}^{1+2\gamma}}
	\,\mathrm d s\mathrm d t
	\leqslant C_{\gamma,T}
	\tau^{2(\gamma-\beta)} \nm{v}_{H^\gamma(0,T)}^2.
\end{align*}
This gives \eqref{eq:I2} and thus proves
\[
\nm{\left( I-P_\tau \right)v}_{H^\beta(0,T)}
\leqslant  C_{\beta,\gamma,T}
\tau^{\gamma-\beta}\nm{v}_{H^1(0,T)},
\]
where $1/2<\gamma<1$.

As the proof of the estimate for $Q_\tau$ is similar, we omit it here and conclude the proof of this lemma.
\end{proof}

\begin{lemma}
	\label{lem:PQ-L-inf}
	If $ 0 \leqslant \gamma <1/2 $ and $ v \in {}_0H^{1+\gamma}(0,T) $, then
	\begin{align}
		\nm{\left( I-P_\tau \right)v}_{L^\infty(0,T)}
		\leqslant {}&\frac{C_\gamma	\tau^{\gamma+1/2}}{1-2\gamma }
		\nm{v}_{{}_0H^{1+\gamma}(0,T)},
		\label{eq:PQ-Linf}
	\end{align}
	where the implicit constant $C_\gamma$ is uniformly bounded whenever $\gamma\to1/2$ or $\gamma\to0$.
\end{lemma}

\begin{proof}		
The case  $\gamma=0$ is standard. Below we consider $0<\gamma<1/2$ and let
$1\leqslant j\leqslant J$ be arbitrary. Observe that
\[
\begin{split}
	{}&	\nm{(I- P_\tau) v}_{L^\infty(I_j)}^2
	\leqslant \tau_j
	\int_{I_j} \snm{v'(t)}^2\mathrm d t
	\leqslant{} \tau_j^{1+2\gamma}
	\int_{I_j} (t-t_{j-1})^{-2\gamma} \snm{v'(t)}^2\mathrm d t.
\end{split}
\]
We extend $v'$ to $\mathbb R\backslash(t_{j-1},T)$ by zero and denote it by $Ev'$.
According to \cite[Lemma 4.3]{li_space-time_2019} and \cite[Lemma 16.3]{Tartar2007}, it holds that
\[
\begin{split}
	{}&	\int_{I_j} (t-t_{j-1})^{-2\gamma} \snm{v'(t)}^2\mathrm d t\\
	\leqslant {}&2\gamma
	\int_{\mathbb R} \int_{\mathbb R}
	\frac{\snm{(Ev')(s)-(Ev')(t)}^2}{\snm{s-t}^{1+2\gamma}}	
	\,\mathrm d s\mathrm d t
	\leqslant{} C\snm{v'}_{H^{\gamma}(t_{j-1},T)}^2,
\end{split}
\]
where $C>0$ is independent of $\gamma,v,j$ and $T$. According to \eqref{eq:norm-1/2}, it is not hard to find $\snm{v'}_{H^{\gamma}(t_{j-1},T)}\leqslant
\snm{v'}_{H^{\gamma}(0,T)}$.
Moreover, using Lemmas~\ref{lem:coer} and \ref{lem:coer1} gives
\[
\begin{split}
	{}&	\snm{v'}_{H^{\gamma}(t_{j-1},T)}^2\leqslant
	\snm{v'}_{H^{\gamma}(0,T)}^2=
	\sec\gamma\pi
	\dual{\D_{0+}^{\gamma}v',\,\D_{T-}^{\gamma}v'}_{(0,T)}\\
	\leqslant{}& \sec^2\gamma\pi\nm{\D_{0+}^{\gamma}v'}_{L^2(0,T)}^2
	\leqslant \frac{C_\gamma}{(1-2\gamma)^2}
	\nm{v}^2_{{}_0H^{1+\gamma}(0,T)},
\end{split}
\]
where the constant $C_\gamma$ is uniformly
bounded whenever $\gamma\to1/2$ or $\gamma\to0$.
This establishes \eqref{eq:PQ-Linf} and completes the proof.
\end{proof}

\begin{lemma}
	\label{lem:foo}
	If $ v \in L^2(0,T;\dot H^1(\Omega)) $ and $ v' \in L^2(0,T;\dot H^{-1}(\Omega)) $,
	then
	\begin{align*}
		\dual{v',V}_{L^2(0,t_j;\dot H^1(\Omega))} &=
		\dual{v(t_j), V_j^{-}}_\Omega -
		\sum_{i=0}^{j-1} \dual{\jmp{V_i}, (Q_\tau P_hv)_i^{+}}_\Omega, \\
		\dual{v',V}_{L^2(0,t_j;\dot H^1(\Omega))} &=
		\sum_{i=0}^{j-1} \dual{\jmp{(P_\tau P_h v)_i},V_i^{+}}_\Omega -
		\dual{v(0),V_0^{+}}_\Omega,
	\end{align*}
	for all $ V \in \mathcal  X_{\tau,h} $ and $ 1 \leqslant j \leqslant J $.
\end{lemma}

\begin{proof}
By \cite[Theorem 3 in \S 5.9.2]{Evansbook}, we have $v\in C([0,T];L^2(\Omega))$ and the integration by parts formula holds
\[
\dual{v',V}_{L^2(t_{j-1},t_{j};\dot H^1(\Omega))} =
\dual{v(t_j), V_j^{-}}_\Omega -
\dual{v(t_{j-1}), V_{j-1}^{+}}_\Omega.
\]
In view of the definitions of $P_h,P_\tau$ and $Q_\tau$, it is not hard to establish the desired results.
\end{proof}

\begin{lemma}
	\label{lem:conv}
	If $ u_0=0$ and $ f \in L^2(0,T;\dot H^{-1}(\Omega)) $, then
	\begin{equation}
		\label{eq:conv_theta}
		\begin{aligned}
			& \nm{U-P_\tau P_hu}_{L^\infty(0,T;L^2(\Omega))} +
			\snm{U-P_\tau P_h u}_{
				H^{(1-\alpha)/2}(0,T;\dot H^1(\Omega))
			} \\
			\leqslant{} & C_{\alpha}
			\snm{(I-P_\tau P_h)u}_{H^{(1-\alpha)/2}(0,T;\dot H^1(\Omega))}.
		\end{aligned}
	\end{equation}
\end{lemma}

\begin{proof}
Set $ \theta = U - P_\tau P_h u $.
By Theorem~\ref{thm:regu_pde}, we have
\begin{equation}\label{eq:uV}
	\dual{u',V}_{L^2(0,T;\dot H^1(\Omega))} +
	\dual{
		\D_{0+}^{1-\alpha}\nabla u,
		\nabla V
	}_{\Omega_T} =
	\dual{f,V}_{L^2(0,T;\dot H^1(\Omega))},
\end{equation}
for all $ V \in \mathcal  X_{\tau,h} $.
We use Lemma~\ref{lem:foo} to rewrite the first term and obtain
\[
\sum_{i=0}^{j-1} \dual{\jmp{(P_\tau P_h u)_i},\theta_i^{+}}_\Omega +
\dual{
	\D_{0+}^{1-\alpha}\nabla u,
	\nabla V
}_{\Omega_T} =
\dual{f,V}_{L^2(0,T;\dot H^1(\Omega))},
\]
which gives the identity
\[
\mathcal A(P_\tau P_hu,V) =
\dual{f,V}_{L^2(0,T;\dot H^1(\Omega))}+
\dual{
	\D_{0+}^{1-\alpha}\nabla (P_\tau P_hu-u),
	\nabla V
}_{\Omega_T}.
\]
This together with \eqref{eq:numer_sol} yields the error equation
\begin{equation}\label{eq:error-eq}
	\mathcal A(\theta,V) =
	\dual{
		\D_{0+}^{1-\alpha}\nabla (u-P_\tau P_hu),
		\nabla V
	}_{\Omega_T},
\end{equation}
for all $ V \in \mathcal  X_{\tau,h} $. Letting $V = \theta\chi_{(0,t_j)}$ and using \eqref{eq:sj} and Lemma~\ref{lem:coer}, we get the estimate
\begin{align*}
	& \nm{\theta^-_j}_{L^2(\Omega)}^2 + \nm{\theta_0^{+}}_{L^2(\Omega)}^2 +
	\snm{\theta}_{H^{(1-\alpha)/2}(0,t_j;\dot H^1(\Omega))}^2 \\
	\leqslant{} & C_{\alpha}
	\snm{u-P_\tau P_hu}_{H^{(1-\alpha)/2}(0,T;\dot H^1(\Omega))}
	\snm{\theta}_{H^{(1-\alpha)/2}(0,t_j;\dot H^1(\Omega))},
\end{align*}
and using Young's inequality with $ \epsilon
$ proves \eqref{eq:conv_theta}.
\end{proof}

\subsection{\bf Proof of  Theorem~\ref{thm:stab}} \label{subsec:5.2} 

Let $ 1 \leqslant j \leqslant J $.
Inserting $ V = U \chi_{(0,t_j)} $ into \eqref{eq:numer_sol} and applying \eqref{eq:sj} yield that
\begin{align*}
	&{}
	\frac12 \left(
	\nm{U^-_j}_{L^2(\Omega)}^2 + \nm{U_0^{+}}_{L^2(\Omega)}^2
	\right) +\sin(\pi\alpha/2)
	\snm{U}_{H^{(1-\alpha)/2}(0,t_j; \dot H^1(\Omega))}^2\\
	\leqslant{}&
	\dual{f,U\chi_{(0,t_j)}}_{H^{(1-\alpha)/2}(0,T; \dot H^1(\Omega))} +
	\dual{u_0,U_0^{+}}_{L^2(\Omega)},
\end{align*}
which further implies
\begin{align*}
	&{}
	\nm{U^-_j}_{L^2(\Omega)}^2 +\sin(\pi\alpha/2)
	\snm{U}_{H^{(1-\alpha)/2}(0,t_j; \dot H^1(\Omega))}^2\\
	\leqslant{}&
	\nm{u_0}_{L^2(\Omega)}^2 +2
	\nm{f}_{(H^{(1-\alpha)/2}(0,T;\dot
		H^{1}(\Omega)) )^*}
	\nm{U\chi_{(0,t_j)}}_{H^{(1-\alpha)/2}(0,T; \dot H^1(\Omega))}.
\end{align*}
Noticing \eqref{eq:norm-1/2} we have
\[
\begin{aligned}
	\nm{U\chi_{(0,t_j)}}_{H^{(1-\alpha)/2}(0,T; \dot H^1(\Omega))}\leqslant{}& \!C_{\alpha,T}
	\snm{U\chi_{(0,t_j)}}_{H^{(1-\alpha)/2}(0,T;\dot H^1(\Omega))} \\
	={}&C_{\alpha,T}
	\snm{U}_{H^{(1-\alpha)/2}(0,t_j;\dot H^1(\Omega))},
\end{aligned}
\]
and invoking Young's inequality with $ \epsilon $, we obtain
\begin{align*}
	{}&	\nm{U_j^{-}}_{L^2(\Omega)}^2 +
	\snm{U}_{H^{(1-\alpha)/2}(0,t_j;\dot H^1(\Omega))}^2\\
	\leqslant{}&C_{\alpha,T}\left(
	\nm{u_0}_{L^2(\Omega)}^2 +
	\nm{f}_{(H^{(1-\alpha)/2}(0,T;\dot
		H^{1}(\Omega)) )^*}^2\right).
\end{align*}
Consequently, it follows that
\[
\begin{aligned}
	{}&	\nm{U}_{L^\infty(0,T;L^2(\Omega))} +
	\snm{U}_{H^{(1-\alpha)/2}(0,T;\dot H^1(\Omega))}\\
	\lesssim{}&
	\nm{u_0}_{L^2(\Omega)} +
	\nm{f}_{(H^{(1-\alpha)/2}(0,T;\dot
		H^{1}(\Omega)) )^*}.\\
\end{aligned}
\]
This proves \eqref{eq:stab} and thus completes the proof.
\subsection{\bf Proof of  Theorem~\ref{thm:conv-f-L2}} \label{subsec:5.3} 

By Theorem~\ref{thm:regu_pde}, Lemmas~\ref{lem:PQ-Hs}, \ref{lem:PQ-L-inf} and \ref{lem:interp}, a routine calculation gives the following estimates
\[
\begin{split}
	\nm{(I\!-\!P_h)u}_{L^{\infty}(0,T;L^2(\Omega))}+
	\snm{(I\!-\!P_h)u}_{H^{(1-\alpha)/2}(0,T;\dot H^1(\Omega))} \lesssim {}&h  	\nm{f}_{L^2(\Omega_T)},\\
	\nm{(I\!-\!P_\tau)u}_{L^{\infty}(0,T;L^2(\Omega)}  \!+\!
	\snm{(I\!-\!P_\tau)P_hu}_{H^{(1-\alpha)/2}(0,T;\dot H^1(\Omega))} \lesssim {}&\tau^{1/2} 	\nm{f}_{L^2(\Omega_T)}.
\end{split}
\]
Therefore, applying Lemma~\ref{lem:conv} yields that
\begin{align*}
	& \nm{u-U}_{L^{\infty}(0,T;L^2(\Omega))} +
	\snm{u-U}_{H^{(1-\alpha)/2}(0,T;\dot H^1(\Omega))} \\
	\leqslant{} &\nm{u-P_\tau P_h u}_{L^{\infty}(0,T;L^2(\Omega))} +
	\nm{U-P_\tau P_h u}_{L^{\infty}(0,T;L^2(\Omega))} \\
	& \quad {} + \snm{u-P_\tau P_h u}_{H^{(1-\alpha)/2}(0,T;\dot H^1(\Omega))} +
	\snm{U-P_\tau P_hu}_{H^{(1-\alpha)/2}(0,T;\dot H^1(\Omega))} \\
	\lesssim {} &\nm{u-P_\tau P_h u}_{L^{\infty}(0,T;L^2(\Omega))} +
	\snm{u-P_\tau P_hu}_{H^{(1-\alpha)/2}(0,T;\dot H^1(\Omega))}\\
	\lesssim{} &
	\big(h + \tau^{1/2}\big)	
	\nm{f}_{L^2(\Omega_T)}.
\end{align*}
This establishes \eqref{eq:conv-f-L2-Hs}.

Next, let us prove \eqref{eq:conv-f-L2-L2} by duality argument.
By Theorem~\ref{thm:dual}, there exists a unique $z\in  \mathcal G$
such that
\[
-\dual{z', v}_{L^2(0,T;\dot H^1(\Omega))} +
\dual{
	\D_{T-}^{1-\alpha}\nabla  z, \nabla v
}_{\Omega_T} =
\dual{u-U, v}_{\Omega_T}
\]
for all $ v \in L^2(0,T;\dot H^1(\Omega)) $. Letting $ v = u-U $ gives
\begin{align*}
	{}& \nm{u-U}^2_{L^2(\Omega_T)} \\
	={}&  - \dual{z',u-U}_{L^2(0,T;\dot H^1(\Omega))} +
	\dual{
		\D_{T-}^{1-\alpha}\nabla  z, \nabla( u-U)
	}_{\Omega_T} \\
	={}&
	\dual{u',z}_{L^2(0,T;\dot H^1(\Omega))} + \dual{z',U}_{L^2(0,T;\dot H^1(\Omega))} +
	\dual{\D_{0+}^{1-\alpha}\nabla(u-U), \nabla z}_{\Omega_T},
\end{align*}
by integration by parts and Lemma~\ref{lem:coer}. Moreover, setting $ Z =Q_\tau P_h
z $ and combining \eqref{eq:numer_sol}, \eqref{eq:uV}, and Lemma~\ref{lem:foo} yield that
\[
\begin{aligned}
	{}&\dual{u',Z}_{L^2(0,T;\dot H^1(\Omega))} +
	\dual{\D_{0+}^{1-\alpha}\nabla(u-U),\nabla Z}_{\Omega_T} \\
	={}&
	\sum_{j=0}^{J-1} \dual{\jmp{U_j},Z_j^{+}}_\Omega=-
	\dual{z',U}_{\Omega_T}.
\end{aligned}
\]
Consequently, we obtain
\begin{align}
	& \nm{u-U}_{L^2(\Omega_T)}^2 \notag\\
	= &
	\dual{u',z-Z}_{L^2(0,T;\dot H^1(\Omega))} +
	\dual{\D_{0+}^{1-\alpha}\nabla(u-U), \nabla(z-Z)}_{\Omega_T} \notag\\
	\leqslant{} &
	\nm{u'}_{L^2(\Omega_T)}
	\nm{z-Z}_{L^2(\Omega_T)}\notag\\
	& \quad {} + \snm{u-U}_{H^{(1-\alpha)/2}(0,T;\dot H^1(\Omega))}
	\snm{z-Z}_{H^{(1-\alpha)/2}(0,T;\dot H^1(\Omega))},\label{eq:dual-u-U}
\end{align}
by the Cauchy\textendash Schwartz inequality and Lemma~\ref{lem:coer}.
In view of Theorem~\ref{thm:dual}, Lemmas~\ref{lem:PQ-Hs} and \ref{lem:interp}, a direct computation implies
\begin{align*}
	\nm{z-Z}_{L^2(\Omega_T)} &
	\lesssim \left( h^{2} + \tau \right)
	\nm{u-U}_{L^2(\Omega_T)}, \\
	\snm{z-Z}_{H^{(1-\alpha)/2}(0,T;\dot H^1(\Omega))} &
	\lesssim \left( h + \tau^{1/2} \right)
	\nm{u-U}_{L^2(\Omega_T)},
\end{align*}
and by \eqref{eq:conv-f-L2-Hs} and Theorem~\ref{thm:regu_pde}, it follows that
\begin{equation*}
	\begin{aligned}
		\nm{u-U}_{L^2(\Omega_T)}
		\lesssim{} &
		\left( h^{2}+\tau \right)
		\nm{f}_{L^2(\Omega_T)}  +		\big( h+\tau^{1/2} \big)
		\big( h + \tau^{1/2} \big)
		\nm{f}_{L^2(\Omega_T)}\\
		\lesssim{}&\left( h^{2}+\tau \right)
		\nm{f}_{L^2(\Omega_T)}.
	\end{aligned}
\end{equation*}
This proves \eqref{eq:conv-f-L2-L2} and concludes the proof of Theorem~\ref{thm:conv-f-L2}.

\subsection{\bf Proof of  Theorem~\ref{thm:conv-u0-L2}} \label{subsec:5.4} 

According to Theorems~\ref{thm:regu-fpde-u0} and \ref{thm:stab} we have the stability estimate
\[
\begin{aligned}
	\nm{U}_{L^\infty(0,T;L^2(\Omega))}
	+\nm{u'}_{(H^{(1-\alpha)/2}(0,T;\dot H^{1}(\Omega))^*}	{}&\\
	\qquad+ \nm{u-U}_{H^{(1-\alpha)/2}(0,T;\dot H^1(\Omega))}
	{}&\lesssim	\nm{u_0}_{L^2(\Omega)}.
\end{aligned}
\]
Recalling the proof of \eqref{eq:conv-f-L2-L2}, we claim that there exists a unique
\[
w\in {}^0\!H^{1}(0,T;L^2(\Omega))\cap
{}^0\!H^{1-\alpha}(0,T;\dot H^2(\Omega))
\]
satisfying
\[
\snm{w-W}_{H^{(1-\alpha)/2}(0,T;\dot H^1(\Omega))}
\lesssim \big(h + \tau^{1/2}\big)
\nm{u-U}_{L^2(\Omega_T)}
\]
and
\[
\begin{split}
	\nm{u\!-\!U}_{L^2(\Omega_T)}^2 \!
	\leqslant&
	\nm{u'}_{(H^{(1-\alpha)/2}(0,T;\dot H^{1}(\Omega))^*}
	\nm{w\!-\!W}_{H^{(1-\alpha)/2}(0,T;\dot H^{1}(\Omega)} \\
	& \quad {} + \snm{u\!-\!U}_{H^{(1-\alpha)/2}(0,T;\dot H^1(\Omega))}
	\!	\snm{w\!-\!W}_{H^{(1-\alpha)/2}(0,T;\dot H^1(\Omega))},
\end{split}
\]
where $W = Q_\tau P_h w$. Consequently, we obtain \eqref{eq:conv-u0-L2-L2} and complete the proof of Theorem~\ref{thm:conv-u0-L2}.

\section{Proof of Theorem~\ref{thm:conv-Linf}} \label{sec:6}

\setcounter{section}{6} \setcounter{equation}{0}

In this section, based on some estimates established in  \cite{McLean2015Time}, we aim to prove  Theorem~\ref{thm:conv-Linf} under further assumption that temporal grid is uniform, i.e., $\tau_j =\tau = T/J$ for all $1\leqslant j\leqslant J$.	

\subsection{\bf Auxiliary results} \label{subsec:6.1} 
Assuming that $ y_0 \in \mathbb R $ and $ \lambda>0$ is a constant, we set $Y_0=y_0$ and define
a sequence $ \{Y_k\}_{k=1}^\infty $ as follows
\begin{equation}
	\label{eq:Y_k}
	\mu \bigg(
	\sum_{j=1}^k \big(
	b_{k-j+2} - 2b_{k-j+1} + b_{k-j}
	\big) Y_j + b_1 Y_{k+1}
	\bigg) + Y_{k+1} - Y_k = 0,
\end{equation}
where $ \mu := \lambda \tau^\alpha $ and $ b_j := j^\alpha/\Gamma(1+\alpha)
$ for all $ j \in \mathbb N $.

\begin{lemma}
	\label{lem:ode-jump}
	The sequence $\{Y_k\}_{k=1}^\infty$ defined by \eqref{eq:Y_k} satisfies that
	\begin{align}
		\snm{Y_k}+		C_\alpha k\snm{Y_{k+1} - Y_k} \leqslant {}&
		\snm{y_0}\quad\forall\,k\geqslant 1,
		\label{eq:ode-jump}
	\end{align}
	with some positive constant $C_\alpha>0$.
\end{lemma}

To prove the above lemma, we shall introduce an auxiliary function
\cite{McLean2015Time}
\[
\psi(z) := \frac{(e^z-1)}{2\pi {\bf i}}
\int_{-\infty}^{0+} \frac{w^{-1-\alpha}}{e^{z-w}-1} \, \mathrm{d}w,
\quad  z \in \mathbb C \setminus (-\infty,0],
\]
where ${\bf i}$ denotes the imaginary unit and $ \int_{-\infty}^{0+} $ means the integration through the Hankel contour, i.e., a smooth and
non-self-intersecting path enclosing the negative real axis and orienting
counterclockwise, $ 0 $ and $ z+2k\pi {\bf i} $, $ k \in \mathbb Z$, lie on the different sides of this path.

By \cite[Lemma 1]{McLean2015Time}, we have the expression
\[
\psi(z) = \frac{1}{2\pi {\bf i}}
\int_{-\infty}^{0+} \frac{w^{-\alpha}}{1-e^{w-z}}\cdot\frac{e^w-1}{w} \, \mathrm{d}w,
\]
and
$1+\nu\psi(z)\neq0$ for all $\nu>0$ and $z \in \mathbb C \setminus (-\infty,0]$.
Since the integrand behaves like $O(w^{-\alpha})$ as $w\to0$ and $O(w^{-\alpha-1})$ as $w\to-\infty$, we can shrink the Hankel contour on to $\{se^{+\pi{\bf i}}:0<s<\infty\}$ and $\{se^{-\pi{\bf i}}:0<s<\infty\}$. This allows us to define
\[
\psi(se^{+ \pi{\bf i}}) : = 		\lim_{\substack{z \to -s\\ \mathrm{Im}\, z > 0}}		\psi(z),\quad
\psi(se^{- \pi{\bf i}}) : = 		\lim_{\substack{z \to -s\\ \mathrm{Im}\,z < 0}}		\psi(z),
\]
where $ se^{+\pi{\bf i}} $ and $ se^{-\pi{\bf i}} $ are identified as two different numbers. 

\medskip

{\bf Proof of Lemma~\ref{lem:ode-jump}.}
We first prove $\snm{Y_k}\leqslant \snm{y_0}$. Note that
\begin{equation}\label{eq:Yk}
	Y_{k+1} = \frac{Y_{k}}{1+\mu b_1}- \frac{\mu}{1+\mu b_1}\sum_{j=1}^{k}
	\big(b_{k-j+2} - 2b_{k-j+1} + b_{k-j}
	\big) Y_{j}.
\end{equation}
Observing the fact
\[
0<(j+1)^\alpha-j^\alpha< j^\alpha-(j-1)^\alpha\leqslant1,\quad\,j>0,
\]
we have the estimate
\[
\begin{split}
	\sum_{j=1}^{k}\snm{b_{k-j+2} - 2b_{k-j+1} + b_{k-j}}
	={}&\!\!\sum_{j=1}^{k}\big(b_{k-j+1} -b_{k-j+2}+b_{k-j+1} - b_{k-j}\big)\\
	={}&b_1+b_k-b_{k+1}\leqslant b_1,
\end{split}
\]
for all $k\in\mathbb N$. Since $Y_0=y_0$, by \eqref{eq:Yk} and the above estimate, a standard deduction gives the stability result $\snm{Y_k}\leqslant \snm{y_0}$.

Then let us establish \eqref{eq:ode-jump}. By \cite[Theorem 3]{McLean2015Time} we have the relation
\begin{align*}
	Y_{k+1} - Y_k = \frac{y_0}{2\pi {\bf i}}
	\int_{-\infty}^{0+} \frac{e^{kz}}{1+\mu \psi(z)}
	\, \mathrm{d}z,
\end{align*}
and using the same technique as that used to derive
\cite[(37)]{McLean2015Time} yields
\begin{align*}
	Y_{k+1} - Y_k &= \frac{y_0}{2\pi {\bf i}}
	\int_0^\infty e^{-ks} \left(
	\frac1{1+\mu\psi(se^{-\pi{\bf i}})} -
	\frac1{1+\mu\psi(se^{+\pi{\bf i}})}
	\right) \, \mathrm{d}s \\
	&=  \frac{y_0 \sin\alpha\pi}{\pi}
	\int_0^\infty \frac{\mu e^{-ks} s^{-\alpha}}{\snm{1+\mu\psi(se^{+\pi{\bf i}})}^2}
	\frac{e^{-s}-1}{s} \, \mathrm{d}s.
\end{align*}
By \cite[Lemma 9]{McLean2015Time}, there exists a positive constant $C_\alpha>0$ such that
\[
\frac{\snm{1+\mu\psi(se^{+\pi{\bf i}})}^{2}}
{1+\mu^2s^{-2\alpha}}
\geqslant C_{\alpha},
\quad 0<s<\infty,
\]
from which we obtain
\begin{equation}\label{eq:diff-Yk}
	\snm{Y_{k+1} - Y_k} \leqslant {}\frac{\snm{y_0}\sin\alpha\pi}{\pi C_\alpha}
	\int_0^\infty \frac{\mu e^{-ks} }
	{s^{\alpha}+\mu^2s^{-\alpha}}\cdot\frac{1-e^{-s}}{s} \, \mathrm{d}s.
\end{equation}
Observing the inequality $(1-e^{-s})/s\leqslant 1$ for all $0<s<\infty$,
we estimate the integral as follows
\[
\begin{split}
	\int_0^\infty \frac{\mu e^{-ks}  }
	{s^{\alpha}+\mu^2 s^{-\alpha}} \cdot\frac{1-e^{-s}}{s} \, \mathrm{d}s
	\leqslant{}&
	\int_0^\infty \frac{\mu e^{-ks} }
	{s^{\alpha}+\mu^2 s^{-\alpha}}\, \mathrm{d}s
	\leqslant{}\int_0^\infty \frac{\mu e^{-ks} }
	{2\mu}\, \mathrm{d}s =\frac{1}{2k}.
\end{split}
\]
Putting this back to \eqref{eq:diff-Yk} proves \eqref{eq:ode-jump} and thus completes the proof of this lemma.
\hfill\ensuremath{\blacksquare}

\medskip 

Define the discrete Laplacian operator $-\Delta_h: S_h\to S_h$ as follows
\[
\dual{-\Delta_hw_h,v_h}_{\Omega}=
\dual{\nabla w_h,\nabla v_h}_{ \Omega}
\quad\forall\,v_h\in S_h,
\]
for all $w_h\in S_h$. It is well-known that $-\Delta_h$ admits an
orthonormal basis $\{\phi_h^n:1
\leqslant n\leqslant  \snm{ S_h}\}$ such that
$-\Delta_h\phi_h^n =\lambda_h^n\phi_h^n$, where $\snm{ S_h}={\rm dim}\, S_h$ and $\{\lambda_h^n:1
\leqslant n\leqslant  \snm{ S_h}\}$ is a non-decreasing positive sequence.
Furthermore, we introduce the average interpolation operator $\Pi_\tau:L^{1}(0,T)\to \mathcal X_{\tau}$ by that
\[
(\Pi_\tau v)|_{I_j} = \frac{1}{\tau_j}\int_{I_j}v(t)\,\mathrm d t
\quad\forall\,1\leqslant j\leqslant J,
\]
where $\mathcal X_\tau$ is given by
\[
\mathcal X_{\tau} := \left\{
v_\tau \in L^2(0,T):\
v_\tau|_{I_j} \in P_0(I_j)
\quad \forall\,1\leqslant j\leqslant J
\right\}.
\]
For the operator $\Pi_\tau$, we have the commutativity
\begin{equation}\label{eq:comm}
	\int_0^Tw(t)(\Pi_\tau v )(t)\,\mathrm d t
	= \int_0^Tv(t)(\Pi_\tau w )(t)\,\mathrm d t,
\end{equation}
for all $w,v\in L^1(0,T)$, and the following estimate is standard
\cite{Ciarlet2002}: if $ 0 \leqslant \beta < 1/2 $ and $ \beta \leqslant \gamma <3/2 $, then
\begin{equation}\label{eq:est-Pi}
	\nm{(I-\Pi_\tau)v}_{H^\beta(0,T)} \lesssim
	\tau^{\gamma-\beta} \nm{v}_{H^\gamma(0,T)}
	\quad \forall \,v \in H^\gamma(0,T).
\end{equation}

Given any $w_h\in S_h$, define $ W \in \mathcal X_{\tau,h} $ such that
\begin{equation}
	\label{eq:W}
	\mathcal A(W,V)
	= \dual{w_h,V_0^+}_{\Omega}
	\quad \forall \,V \in \mathcal  X_{\tau,h}.
\end{equation}
The well-posedness of the above problem follows directly from Theorem~\ref{thm:stab},
and thanks to Lemma~\ref{lem:ode-jump}, we have a stability result which is crucial to our error analysis.

\begin{lemma}
	For any $w_h\in S_h$, the unique solution $W\in \mathcal X_{\tau,h}$ to \eqref{eq:W} satisfies that
	\begin{align}
		\nm{W}_{L^{\infty}(0,T;L^2(\Omega))} +
		\snm{W}_{H^{(1-\alpha)/2}(0,T;\dot H^1(\Omega))} &
		\leqslant C_\alpha \nm{w_h}_{L^2(\Omega)},
		\label{eq:est-W-1} \\
		\sum_{j=1}^{J-1}\nm{\jmp{W_j}}_{L^2(\Omega)}\!\!+\!\!
		\nm{\Pi_\tau\!\D_{0+}^{1-\alpha}\!\!\Delta_hW}_{L^{1}(0,T;L^2(\Omega))}&\!
		\!		\leqslant \!C_{\alpha}\! \big( 1\!\!+\!\ln  J \big)\! \nm{w_h}_{L^2(\Omega)}.
		\label{eq:est-W-2}
	\end{align}
\end{lemma}

\begin{proof}
Let us first prove \eqref{eq:est-W-1}.
By \eqref{eq:sj}, for any $1\leqslant j\leqslant J$, inserting $ V = W\chi_{(0,t_j)} $ into \eqref{eq:W} implies
\[\nm{W_j^{-}}_{L^2(\Omega)}+
\snm{W}_{H^{(1-\alpha)/2}(0,t_j;\dot H^1(\Omega))}
\leqslant C_\alpha \nm{w_h}_{L^2(\Omega)},\quad
1\leqslant j\leqslant J.
\]
Therefore, \eqref{eq:est-W-1} is obtained directly from the above estimate.

Then let us prove \eqref{eq:est-W-2}. It is not hard to find that the solution to \eqref{eq:W} has explicit expression $W = \sum_{n=1}^{\snm{ S_h}}w_n\phi_n$. Here $w_n\in \mathcal X_\tau$ and $w_n|_{I_j}=Y^n_j$ for all $1\leqslant j\leqslant J$, where $\{Y^n_j\}_{j=1}^J$ satisfies \eqref{eq:Y_k} with $Y_0^n =\dual{w_h,\phi_n}_{\Omega}$ and $\mu = \lambda_h^n\tau^\alpha$.
By Lemma~\ref{lem:ode-jump}, we have
\[
\snm{Y_{j}^n-Y_{j-1}^n}\leqslant \frac{C_\alpha}{j}\snm{Y_0^n}.
\]
Hence, it follows that
\[
\nm{\jmp{W_j}}_{L^2(\Omega)}\leqslant
\frac{C_\alpha}{j}\nm{w_h}_{L^2(\Omega)},
\quad1\leqslant j\leqslant J,
\]
which yields the estimate
\begin{equation}\label{eq:jmp}
	\sum_{j=1}^{J-1}
	\nm{\jmp{W_j}}_{L^2(\Omega)}
	\leqslant{}
	C_\alpha \sum_{j=1}^{J-1} j^{-1} 	\nm{w_h}_{L^2(\Omega)}
	\leqslant C_\alpha\big( 1+\ln J \big)  \nm{w_h}_{L^2(\Omega)}.
\end{equation}
Define $V\in\mathcal X_{\tau,h}$ by that
\[
V|_{I_{1}} = W_0^+-w_h ,\quad V|_{I_{j}} = \jmp{W_{j-1}},\quad
2\leqslant j\leqslant J,
\]
then plugging  $V\chi_{(I_j)}$ into \eqref{eq:W} leads to
\[
\tau_1(\Pi_\tau\D_{0+}^{1-\alpha}\Delta_h W)\big|_{I_1} =  W_0^+-w_h,
\quad
\tau_j(\Pi_\tau\D_{0+}^{1-\alpha}\Delta_h W)\big|_{I_j} =  \jmp{W_{j-1}},
\]
for all $2\leqslant j\leqslant J$.
Hence, by \eqref{eq:est-W-1} and \eqref{eq:jmp}, we obtain that
\[
\begin{split}
	\nm{\Pi_\tau\D_{0+}^{1-\alpha}\Delta_hW}_{L^{1}(0,T;L^2(\Omega))}
	={}&\nm{W_0^+\!-\!w_h}_{L^2(\Omega)} +
	\sum_{j=1}^{J-1}\nm{\jmp{W_j}}_{L^2(\Omega)} \\
	\leqslant{}& C_\alpha\big( 1+\ln J \big)  \nm{w_h}_{L^2(\Omega)},
\end{split}
\]
which together with \eqref{eq:jmp} establishes \eqref{eq:est-W-2} and concludes the proof of this lemma.
\end{proof} 

Symmetrically, we have the following lemma. As the proof is similar, we omit it here.

\begin{lemma}
	\label{lem:dual-jump}
	Given any $w_h\in S_h$, if
	$ W \in \mathcal X_{\tau,h} $ satisfies that
	\[
	\mathcal A(V,W)	= \dual{w_h,V_J^-}_{\Omega}
	\quad \forall\,V \in \mathcal X_{\tau,h},
	\]
	then the following estimates hold:
	\begin{align*}
		\nm{W}_{L^{\infty}(0,T;L^2(\Omega))} +
		\snm{W}_{H^{(1-\alpha)/2}(0,T;\dot H^1(\Omega))} &
		\leqslant C_\alpha \nm{w_h}_{L^2(\Omega)},\\
		\sum_{j=1}^{J-1}\nm{\jmp{W_j}}_{L^2(\Omega)}\!+\!
		\nm{\Pi_\tau\D_{T-}^{1-\alpha}\Delta_hW}_{L^{1}(0,T;L^2(\Omega))}&
		\!		\leqslant \!C_{\alpha} \big( 1\!+\!\ln J \big) \nm{w_h}_{L^2(\Omega)}.
	\end{align*}
\end{lemma}
\subsection{\bf Main proof} \label{subsec:6.2} 
We now arrive at a position for proving Theorem~\ref{thm:conv-Linf}.
\vskip0.2cm
We first prove \eqref{eq:conv-LinfL2-f-L2}. To do this, we shall establish the estimate
\begin{equation}\label{eq:todo}
	\begin{split}
		\nm{\theta}_{L^\infty(0,T;L^2(\Omega))}\lesssim{}&
		\nm{(I-\Pi_\tau)u}_{H^{(1-\alpha)/2}(0,T;\dot H^1(\Omega))}\\
		{}&\qquad+
		(1+\snm{\ln\tau} )
		\nm{R_hu-P_\tau P_hu}_{L^{\infty}(0,T;L^2(\Omega))},
	\end{split}
\end{equation}
where $ \theta := U - P_\tau P_h u $. Let $W\in \mathcal X_{\tau,h}$ be defined as follows
\[
\mathcal A(V,W)	= \dual{\theta_J^{-},V_J^-}_{\Omega}
\quad \forall\,V \in \mathcal X_{\tau,h}.
\]
Plugging $V = \theta$ into the above equation and observing the error equation \eqref{eq:error-eq}, we obtain the identity
\begin{align*}
	\nm{\theta_J^{-}}_{L^2(\Omega)}^2
	={}&\mathcal A(\theta,W)	=
	\dual{
		\D_{0+}^{1-\alpha}\nabla(u-P_\tau P_h u), \nabla W
	}_{\Omega_T} \\
	={}&	\dual{
		\D_{0+}^{1-\alpha}\nabla(R_hu-P_\tau P_h u), \nabla W
	}_{\Omega_T} \\
	={} &\dual{
		P_\tau P_hu-R_hu,\, \D_{T-}^{1-\alpha}\Delta_h W
	}_{\Omega_T}.
\end{align*}
Thanks to the commutativity \eqref{eq:comm}, we have
\[ 
\begin{split}
	{}&\nm{\theta_J^{-}}_{L^2(\Omega)}^2 \\
	= &
	\dual{P_\tau P_hu\!-\!R_hu,\,(I\!-\!\Pi_\tau) \D_{T-}^{1-\alpha}\Delta_h W
	}_{\Omega_T}\!\!+\!
	\dual{	P_\tau P_hu\!-\!R_hu,\, \Pi_\tau\!
		\D_{T-}^{1-\alpha}\Delta_h W}_{\Omega_T}\\
	= &-
	\dual{R_hu,\,(I-\Pi_\tau) \D_{T-}^{1-\alpha}\Delta_h W
	}_{\Omega_T}+
	\dual{	P_\tau P_hu-R_hu,\, \Pi_\tau
		\D_{T-}^{1-\alpha}\Delta_h W}_{\Omega_T}\\
	=&\dual{(I-\Pi_\tau)R_hu,\, \D_{T-}^{1-\alpha}(-\Delta_h W)
	}_{\Omega_T}+
	\dual{	P_\tau P_hu-R_hu,\, \Pi_\tau
		\D_{T-}^{1-\alpha}\Delta_h W}_{\Omega_T}.
\end{split}
\] 
Applying Lemmas~\ref{lem:dual-jump} and \ref{lem:coer} yields
\[
\begin{split}
	{}&\dual{	P_\tau P_hu-R_hu,\, \Pi_\tau
		\D_{T-}^{1-\alpha}\Delta_h W}_{\Omega_T}\\
	\lesssim {}&
	\nm{\Pi_\tau\D_{T-}^{1-\alpha}\Delta_h W}_{L^{1}(0,T;L^2(\Omega))}
	\nm{R_hu-P_\tau P_hu}_{L^{\infty}(0,T;L^2(\Omega))}\\
	\lesssim{} &(1+\snm{\ln\tau} )\nm{\theta_J^-}_{L^2(\Omega)}	\nm{R_hu-P_\tau P_hu}_{L^{\infty}(0,T;L^2(\Omega))}
\end{split}
\]
and
\[
\begin{split}
	{}&\dual{(I- \Pi_\tau)R_hu,\, \D_{T-}^{1-\alpha}(-\Delta_h W)}_{\Omega_T}\\
	={}& \dual{(I- \Pi_\tau)\nabla R_hu,\, \D_{T-}^{1-\alpha}\nabla W}_{\Omega_T}\\
	\leqslant {}&\snm{W}_{H^{(1-\alpha)/2}(0,T;\dot H^1(\Omega))}
	\snm{(I-\Pi_\tau)R_hu}_{H^{(1-\alpha)/2}(0,T;\dot H^1(\Omega))}\\
	\lesssim {}&\nm{\theta_J^-}_{L^2(\Omega)}
	\nm{(I-\Pi_\tau)u}_{H^{(1-\alpha)/2}(0,T;\dot H^1(\Omega))}.
\end{split}
\]
Combining the above two estimates gives
\[
\begin{aligned}
	\nm{\theta_J^{-}}_{L^2(\Omega)}\lesssim{}&
	\nm{(I-\Pi_\tau)u}_{H^{(1-\alpha)/2}(0,T;\dot H^1(\Omega))}\\
	{}&\quad+
	(1+\snm{\ln\tau} )
	\nm{R_hu-P_\tau P_hu}_{L^{\infty}(0,T;L^2(\Omega))},\\
\end{aligned}
\]
and similarly one can prove this for $\theta_j^{-}$ with $1\leqslant j<J$. Therefore, the estimate \eqref{eq:todo} follows immediately.

We then estimate the right hand side terms in \eqref{eq:todo}.
By \eqref{eq:est-Pi}, Theorem~\ref{thm:regu_pde}, Lemmas~\ref{lem:PQ-L-inf} and \ref{lem:interp}, we have that
\begin{equation}\label{eq:est-1}
	\begin{aligned}
		\nm{(I\!-\!P_\tau)u}_{L^\infty(0,T;L^2(\Omega))}\!+\!
		\nm{(I\!-\!\Pi_\tau)u}_{H^{\frac{1-\alpha}2}(0,T;\dot H^1(\Omega))}
		\!	\lesssim &\tau^{1/2}\!\nm{f}_{L^2(\Omega_T)},
	\end{aligned}
\end{equation}
and applying Theorem~\ref{thm:regu_pde} and Lemmas~\ref{lem:interp} again implies
\begin{equation}\label{eq:est-2}
	\begin{split}
		{}&\nm{(I-R_h)u}_{L^\infty(0,T;L^2(\Omega))}+
		\nm{(I-P_h)u}_{L^\infty(0,T;L^2(\Omega))}\\
		\lesssim{}&
		\left\{
		\begin{aligned}
			&h^{\min\{2,1/\alpha\}}	\nm{f}_{L^2(\Omega_T)}
			&&\text{ if } \alpha\neq 1/2,\\
			&\frac{1}{\sqrt{\epsilon}}h^{2-\epsilon}	\nm{f}_{L^2(\Omega_T)}
			&&\text{ if } \alpha= 1/2,
		\end{aligned}
		\right.
	\end{split}
\end{equation}
where $0<\epsilon\leqslant 1$. For $\alpha=1/2$, we choose $\epsilon = 1/(1+\snm{\ln h})$ to obtain
\begin{equation}\label{eq:est-3}
	\begin{aligned}
		{}&	\nm{(I-R_h)u}_{L^\infty(0,T;L^2(\Omega))}+
		\nm{(I-P_h)u}_{L^\infty(0,T;L^2(\Omega))}\\
		\lesssim{}&
		\big(1+\sqrt{\snm{\ln h}}\big)h^2\nm{f}_{L^2(\Omega_T)}.
	\end{aligned}
\end{equation}
Plugging \eqref{eq:est-1}\textendash\eqref{eq:est-3} into \eqref{eq:todo} and using the following estimate
\[
\begin{split}
	{}&\nm{R_hu-P_\tau P_hu}_{L^{\infty}(0,T;L^2(\Omega))}\\
	\leqslant {}&
	\nm{(I-R_h)u}_{L^{\infty}(0,T;L^2(\Omega))}+
	\nm{(I-P_h)u}_{L^{\infty}(0,T;L^2(\Omega))}\\
	{}&\quad	+\nm{(I-P_\tau)P_h u}_{L^{\infty}(0,T;L^2(\Omega))}\\
	\leqslant {}&
	\nm{(I-R_h)u}_{L^{\infty}(0,T;L^2(\Omega))}+
	\nm{(I-P_h)u}_{L^{\infty}(0,T;L^2(\Omega))}\\
	{}&\quad	+\nm{(I-P_\tau) u}_{L^{\infty}(0,T;L^2(\Omega))},
\end{split}
\]
we find that
\begin{equation}\label{eq:theta}
	\nm{\theta}_{L^\infty(0,T;L^2(\Omega))}
	\lesssim{}(1+\snm{\ln \tau})
	\left(\tau^{1/2} + \epsilon_h
	h^{\min\{2,1/\alpha\}}\right)
	\nm{f}_{L^2(\Omega_T)},
\end{equation}
where $\epsilon_h = 1$ if $\alpha\neq1/2$ and $\epsilon_h = 1+\sqrt{\snm{\ln h}}$ if $\alpha=1/2$.
Consequently, \eqref{eq:conv-LinfL2-f-L2} follows from \eqref{eq:est-1}\textendash\eqref{eq:theta} and the estimate
\[
\begin{aligned}
	{}&\nm{u-U}_{L^\infty(0,T;L^2(\Omega))}
	\leqslant \nm{\theta}_{L^\infty(0,T;L^2(\Omega))}\!+\!
	\nm{u-P_\tau P_hu}_{L^\infty(0,T;L^2(\Omega))}\\
	\leqslant &\nm{\theta}_{L^\infty(0,T;L^2(\Omega))}\!+\!
	\nm{(I-P_h)u}_{L^\infty(0,T;L^2(\Omega))}\!+\!
	\nm{(I-P_\tau) u}_{L^\infty(0,T;L^2(\Omega))}.
\end{aligned}
\]

In view of \eqref{eq:est-Pi} and Lemma~\ref{lem:PQ-L-inf}, we can establish the estimate \eqref{eq:conv-LinfL2-f-Hs} similarly. This ends the proof of Theorem~\ref{thm:conv-Linf}.

\section{Numerical experiments} \label{sec:7}
\setcounter{section}{7}

In this section, we present several numerical experiments to verify the theoretical results with $ T=1$ and $\Omega = (0,1) $. We will use uniform
grids both in time and space and introduce the following notations:
\[
\begin{aligned}
	\mathcal E_1 & := \nm{\widehat u-U}_{L^2(\Omega_T)}, \\
	\mathcal E_2 & := \nm{\widehat u-U}_{L^\infty(0,T;L^2(\Omega))},	\\
	\mathcal E_3 & :=  \sqrt{\dual{\D_{0+}^{1-\alpha} (\nabla\widehat{u} - \nabla U),\nabla\widehat{u} - \nabla U}_{\Omega_T}},
\end{aligned}
\]
where the reference solution $ \widehat u $ is the numerical solution with
respect to $h = 2^{-10}$ and $\tau = 2^{-15}$. Note that, by Lemma~\ref{lem:coer},
\[
\mathcal E_3\sim	\nm{\D_{0+}^{(1-\alpha)/2} (\nabla\widehat{u} - \nabla U)}_{L^2(\Omega_T)}
\sim\nm{\widehat{u} - U}_{H^{\frac{1-\alpha}{2}}(0,T;\dot H^1(\Omega))}.
\]

With uniform temporal grids, the DG scheme \eqref{eq:numer_sol} results in a block triangular Toeplitz-like with tri-diagonal block system, and we can adopt the fast direct method proposed in \cite{ke_fast_2015} to solve it efficiently with quasi-optimal complexity $O((\tau h)^{-1}|\ln\tau|^2)$. Moreover, $\mathcal E_3$ can be computed via fast Fourier transform.

\medskip

{\bf{Experiment 1.}}
Consider
\begin{alignat*}{2}
	u_0(x) & := 0, & \quad & x\in\Omega, \\
	f(x,t) & := x^{-0.49} t^{-0.49}, & \quad & (x,t)\in\Omega_T.
\end{alignat*}
To test the accuracy of spatial discretization, we fix temporal step size
$\tau = 2^{-15}$. Since $f\in L^2(\Omega_T)$, according to Theorems~\ref{thm:conv-f-L2} and \ref{thm:conv-Linf}, we have $\mathcal E_1= O(h^{2}),\,
\mathcal E_2= O\big(h^{\min\{2,1/\alpha\}}\big)$ and $\mathcal E_3= O(h)$.
These coincide with the numerical results in Table~\ref{tab:subdf-ex1-space}.
\begin{table}[H] 
	\small\setlength{\tabcolsep}{3pt}
	\begin{center}
		\begin{tabular}{cccccccc}
			\hline\noalign{\smallskip}
			&		 $h$ & $ \mathcal E_1$ & Order & $ \mathcal E_2 $ &
			Order & $ \mathcal E_3 $ & Order\\
			\noalign{\smallskip}\hline\noalign{\smallskip}
			\multirow{4}{*}{$\alpha=0.8$}
			&        $2^{-2}$           & 2.00e-02 & --  & 3.67e-02 & --  & 2.87e-01 & -- \\
			&        $2^{-3}$           & 5.53e-03 & 1.85  & 1.48e-02 & 1.31  & 1.59e-01 & 0.85 \\
			&        $2^{-4}$          & 1.50e-03 & 1.88  & 5.95e-03 & 1.31  & 8.68e-02 & 0.87 \\
			&        $2^{-5}$          & 4.04e-04 & 1.90  & 2.42e-03 & 1.30  & 4.67e-02 & 0.89 \\
			\noalign{\smallskip}\hline\noalign{\smallskip}
			\multirow{4}{*}{$\alpha=0.2$}
			&        $2^{-4}$          & 1.13e-03 &--  & 1.45e-03 & --  & 6.05e-02 & --\\
			&        $2^{-5}$          & 3.02e-04 & 1.90  & 3.88e-04 & 1.91  & 3.21e-02 & 0.91 \\
			&        $2^{-6}$          & 7.99e-05 & 1.92  & 1.02e-04 & 1.92  & 1.69e-02 & 0.93 \\
			&        $2^{-7}$         & 2.08e-05 & 1.94  & 2.67e-05 & 1.94  & 8.81e-03 & 0.94 \\
			\noalign{\smallskip}\hline
		\end{tabular}
		\smallskip
		\caption{Spatial errors of {\bf Experiment 1} with $ \tau = 2^{-15}$.
		}
		\label{tab:subdf-ex1-space}
	\end{center}
\end{table}	

Next, we consider temporal errors and choose $h = 2^{-10}$. In view of
Table~\ref{tab:subdf-ex1-temp}, we find that $\mathcal E_1= O(\tau),\,
\mathcal E_2= O(\tau^{1/2})$ and $\mathcal E_3= O(\tau^{1/2})$.
Evidently, they match well the estimates given by Theorems~\ref{thm:conv-f-L2} and \ref{thm:conv-Linf}.

\begin{table}[H] 
	\small\setlength{\tabcolsep}{3pt}
	\begin{center}
		\begin{tabular}{cccccccc}
			\hline\noalign{\smallskip}
			&		 $\tau$  & $ \mathcal E_1$ & Order & $ \mathcal E_2 $ &
			Order & $ \mathcal E_3 $ & Order\\
			\noalign{\smallskip}\hline\noalign{\smallskip}
			\multirow{4}{*}{$\alpha=0.7$}
			&        $2^{-9}$           & 1.99e-03 & --  & 9.22e-02 & --  & 2.12e-02 & -- \\
			&       $2^{-10}$          & 1.13e-03 & 0.81  & 6.54e-02 & 0.49  & 1.42e-02 & 0.58 \\
			&        $2^{-11}$         & 6.24e-04 & 0.86  & 4.43e-02 & 0.56  & 9.28e-03 & 0.61 \\
			&        $2^{-12}$           & 3.27e-04 & 0.93  & 2.81e-02 & 0.66  & 5.86e-03 & 0.66 \\
			\noalign{\smallskip}\hline\noalign{\smallskip}
			\multirow{4}{*}{$\alpha=0.3$}
			&       $2^{-9}$         & 6.24e-04 & --  & 3.37e-02 & --  & 2.92e-02 & --\\
			&       $2^{-10}$        & 3.63e-04 & 0.78  & 2.49e-02 & 0.43  & 2.10e-02 & 0.48 \\
			&       $2^{-11}$        & 2.06e-04 & 0.82  & 1.77e-02 & 0.50  & 1.48e-02 & 0.51 \\
			&        $2^{-12}$       & 1.12e-04 & 0.88  & 1.17e-02 & 0.59  & 1.00e-02 & 0.56 \\
			\noalign{\smallskip}\hline
		\end{tabular}
		\smallskip
		\caption{Temporal errors of {\bf Experiment 1} with $h = 2^{-10}$.}
		\label{tab:subdf-ex1-temp}
	\end{center}
\end{table}

\begin{table}[H] 
	\small\setlength{\tabcolsep}{3pt}
	\begin{center}
		\begin{tabular}{ccccccccc}
			\hline\noalign{\smallskip}
			&	\multicolumn{2}{c}{$\alpha=0.9$} &&\multicolumn{2}{c}{$\alpha=0.5$} &&\multicolumn{2}{c}{$\alpha=0.3$}\\
			\noalign{\smallskip}\cline{2-3} \cline{5-6} \cline{8-9} \noalign{\smallskip}
			$h	 $   & $ \mathcal E_2 $ & Order & \phantom{aa} &
			$ \mathcal E_2 $ & Order& \phantom{aa} &$ \mathcal E_2 $ & Order \\
			\noalign{\smallskip}\hline\noalign{\smallskip}
			$2^{-4}$          & 7.10e-04 & --  &       & 5.81e-04 & -- &       & 5.18e-04 & -- \\
			$2^{-5}$         & 1.90e-04 & 1.91  &       & 1.55e-04 & 1.90  &       & 1.39e-04 & 1.90 \\
			$2^{-6}$         & 5.01e-05 & 1.92  &       & 4.11e-05 & 1.92  &       & 3.66e-05 & 1.92 \\
			$2^{-7}$        & 1.30e-05 & 1.94  &       & 1.07e-05 & 1.94  &       & 9.55e-06 & 1.94 \\			
			\noalign{\smallskip}\hline
		\end{tabular}
		\smallskip
		\caption{Spatial errors of {\bf Experiment 2} with $ \tau = 2^{-15}$.}
		\label{tab:subdf-ex2-space}
	\end{center}
\end{table}

\medskip

{\bf{Experiment 2.}}
Consider
\begin{alignat*}{2}
	u_0(x) & := 0, & \quad & x\in\Omega, \\
	f(x,t) & := x^{-0.49} t^{0.01}, & \quad & (x,t)\in\Omega_T.
\end{alignat*}
It is clear that $f\in {}_0H^{1/2}(0,T;L^2(\Omega))$. In Tables~\ref{tab:subdf-ex2-space} and \ref{tab:subdf-ex2-temp}, we observe the optimal convergence order
$ \mathcal E_2= O(\tau+h^2) $, which agrees with Theorem~\ref{thm:conv-Linf}.	
\begin{table}[H] 
	\small\setlength{\tabcolsep}{3pt}
	\begin{center}
		\begin{tabular}{ccccccccc}
			\hline\noalign{\smallskip}
			&	\multicolumn{2}{c}{$\alpha=0.7$} &&\multicolumn{2}{c}{$\alpha=0.4$} &&\multicolumn{2}{c}{$\alpha=0.1$}\\
			\noalign{\smallskip}\cline{2-3} \cline{5-6} \cline{8-9} \noalign{\smallskip}
			$\tau	 $   & $ \mathcal E_2 $ & Order & \phantom{aa} &
			$ \mathcal E_2 $ & Order& \phantom{aa} &$ \mathcal E_2 $ & Order \\
			\noalign{\smallskip}\hline\noalign{\smallskip}
			$2^{-8}$       & 3.18e-04 & --  &       & 2.02e-04 & --  &       & 2.14e-04 & -- \\
			$2^{-9}$          & 1.60e-04 & 1.00  &       & 1.00e-04 & 1.01  &       & 1.04e-04 & 1.04 \\
			$2^{-10}$         & 7.95e-05 & 1.01  &       & 4.97e-05 & 1.01  &       & 5.03e-05 & 1.04 \\
			$2^{-11}$         & 3.92e-05 & 1.02  &       & 2.46e-05 & 1.02  &       & 2.43e-05 & 1.05 \\
			\noalign{\smallskip}\hline
		\end{tabular}
		\smallskip
		\caption{Temporal errors of {\bf Experiment 2} with
			$h = 2^{-10}$.}
		\label{tab:subdf-ex2-temp}
	\end{center}
\end{table}

\begin{table}[H] 
	\small\setlength{\tabcolsep}{3pt}
	\begin{center}
		\begin{tabular}{ccccccccc}
			\hline\noalign{\smallskip}
			&	\multicolumn{2}{c}{$\alpha=0.9$} &&\multicolumn{2}{c}{$\alpha=0.6$} &&\multicolumn{2}{c}{$\alpha=0.3$}\\
			\noalign{\smallskip}\cline{2-3} \cline{5-6} \cline{8-9} \noalign{\smallskip}
			$\tau	 $ & $ \mathcal E_1 $ & Order & \phantom{aa} &
			$ \mathcal E_1 $ & Order& \phantom{aa} &$ \mathcal E_1 $ & Order \\
			\noalign{\smallskip}\hline\noalign{\smallskip}
			$2^{-7}$              & 2.90e-02 & --  &       & 2.18e-02 & --  &       & 1.09e-02 & -- \\
			$2^{-8}$            & 2.00e-02 & 0.54  &       & 1.46e-02 & 0.58  &       & 8.07e-03 & 0.44 \\
			$2^{-9}$            & 1.37e-02 & 0.54  &       & 9.77e-03 & 0.58  &       & 5.82e-03 & 0.47 \\
			$2^{-10}$            & 9.36e-03 & 0.55  &       & 6.53e-03 & 0.58  &       & 4.07e-03 & 0.52 \\
			\noalign{\smallskip}\hline
		\end{tabular}
		\smallskip
		\caption{Temporal errors of {\bf Experiment 3} with
			$h = 2^{-10}$.}
		\label{tab:subdf-ex3-temp}
	\end{center}
\end{table}

\medskip 

{\bf{Experiment 3.}}
In third test, let us verify Theorem~\ref{thm:conv-u0-L2} and take
\begin{alignat*}{2}
	u_0(x) & := x^{-0.49} , & \quad & x\in\Omega,\\
	f(x,t) & := 0, & \quad & (x,t)\in\Omega_T.
\end{alignat*}
The convergence rate $ \mathcal E_1
= O(\tau^{1/2})$ in Table~\ref{tab:subdf-ex3-temp} coincides with Theorem~\ref{thm:conv-u0-L2}.
However, as we mentioned in Remark~\ref{rem:opt-rate}, Theorem~\ref{thm:conv-u0-L2}
only gives suboptimal spatial rate $ \mathcal E_1= O(h)$. The optimal order of spatial discretization should be $ \mathcal E_1= O(h^{\min\{2,1/\alpha\}})$, which can be observed from Table~\ref{tab:subdf-ex3-space}.

\begin{table}[H] 
	\small\setlength{\tabcolsep}{3pt}
	\begin{center}
		\begin{tabular}{ccccccccc}
			\hline\noalign{\smallskip}
			&	\multicolumn{2}{c}{$\alpha=0.8$} &&\multicolumn{2}{c}{$\alpha=0.5$} &&\multicolumn{2}{c}{$\alpha=0.2$}\\
			\noalign{\smallskip}\cline{2-3} \cline{5-6} \cline{8-9} \noalign{\smallskip}
			$h	 $  & $ \mathcal E_1 $ & Order & \phantom{aa} &
			$ \mathcal E_1 $ & Order& \phantom{aa} &$ \mathcal E_1 $ & Order \\
			\noalign{\smallskip}\hline\noalign{\smallskip}
			$2^{-2}$            & 3.37e-02 & --  &       & 1.54e-02 & --  &       & 1.10e-02 & -- \\
			$2^{-3}$          & 1.36e-02 & 1.31  &       & 4.49e-03 & 1.78  &       & 3.03e-03 & 1.86 \\
			$2^{-4}$         & 5.31e-03 & 1.36  &       & 1.27e-03 & 1.82  &       & 8.20e-04 & 1.89 \\
			$2^{-5}$       & 1.90e-03 & 1.48  &       & 3.48e-04 & 1.86  &       & 2.19e-04 & 1.90 \\
			\noalign{\smallskip}\hline
		\end{tabular}
		\smallskip
		\caption{Spatial errors of {\bf Experiment 3} with
			$ \tau = 2^{-15}$.}
		\label{tab:subdf-ex3-space}
	\end{center}
\end{table}	 

\medskip 

{\bf{Experiment 4.}}
Although the rate $O(\tau^{1/2})$ established in Theorem~\ref{thm:conv-f-L2} is optimal with respect to the Sobolev regularity, it can be further  improved via graded grids, provided that the solution possesses some growth estimates like \eqref{eq:regu-assum}.

To the end, let us investigate the performance of the DG scheme \eqref{eq:numer_sol} under graded temporal grid $	t_j =(j/J)^\sigma,\,j = 0,1,\cdots,J$,
with $\sigma>1$.
For simplicity we pay attention to the quantity $\mathcal E_2$, which corresponds to the $L^\infty(0,T;L^2(\Omega))$-norm, and consider three cases:
\begin{itemize}
	\item[\textbullet] {\it Case 1}: $u_0(x) = x^{-0.49},\,f (x,t)= 0$;
	\item[\textbullet] {\it Case 2}: $u_0(x)=0,\,f (x,t)= x^{-0.49}t^{-0.49}$;
	\item[\textbullet] {\it Case 3}: $u_0(x)=0,\,f (x,t)= x^{-0.49}\snm{1-2t}^{-0.49}$.
\end{itemize}
Note that for all cases we have $u_0\in L^2(\Omega)$ and $f\in L^2(\Omega_T)$. According to \cite{McLean2010}, one can obtain growth estimates for the first two cases and the first order accuracy $\mathcal E_2 = O(\tau)$ is maintained with suitable parameter $\sigma>1$; see Tables~\ref{tab:ex4-u0} and \ref{tab:ex4-f-L2-1}.

\begin{table}[H] 
	\small\setlength{\tabcolsep}{3pt}
	\begin{center}
		\begin{tabular}{ccccccccccc}
			\hline\noalign{\smallskip}
			\multicolumn{5}{c}{$\alpha = 0.3$}  && \multicolumn{5}{c}{$\alpha=0.9$}    \\
			\noalign{\smallskip}\cline{1-5} \cline{7-11}\noalign{\smallskip}
			$\sigma$&& $J$ & $\mathcal E_2 $ &{Order}    &
			& 		$\sigma$&&$J$ & $\mathcal E_2 $ &{Order} \\
			\noalign{\smallskip}\hline\noalign{\smallskip}
			\multirow{4}{*}{2}	
			&  & $2^5$  & 9.91e-01 & --   &&		\multirow{4}{*}{1.5}&  & $2^5$  &  1.17e-00 & -- \\
			&  & $2^6$  & 8.59e-01 & 0.21 &  &&& $2^6$  & 9.80e-01 & 0.26 \\
			&  &$2^7$ & 7.26e-01 & 0.24 &  &&& $2^7$ & 7.77e-01 & 0.33 \\
			&  &$2^8$ &  5.98e-01 & 0.28 &  & &&$2^8$ & 5.68e-01 & 0.45 \\
			\noalign{\smallskip}\hline\noalign{\smallskip}
			\multirow{4}{*}{5}
			&  & $2^5$  & 1.09e-00 & --   &  &\multirow{4}{*}{2.5}&& $2^5$  &  5.17e-01 & -- \\
			&  & $2^6$  & 8.73e-01& 0.32 &  &&& $2^6$ & 3.09e-01 & 0.74 \\
			&  &$2^7$ & 6.42e-01 & 0.44 &  &&& $2^7$ & 1.67e-01 & 0.89 \\
			&  & $2^8$ &  4.16e-01 & 0.63 &  &&& $2^8$ & 8.62e-02 & 0.95 \\
			\noalign{\smallskip}\hline\noalign{\smallskip}	
			\multirow{4}{*}{9}
			&  & $2^5$  & 3.81e-01& --   &  &\multirow{4}{*}{4}&& $2^5$  &  2.38e-01 & -- \\
			&  & $2^6$  & 2.03e-01 & 0.91 &  &&& $2^6$ & 1.25e-01 & 0.93 \\
			&  &$2^7$ & 1.02e-01& 0.99 &  &&& $2^7$ & 6.24e-02 & 1.00 \\
			&  & $2^8$ &  4.98e-02 & 1.03 &  &&& $2^8$ & 3.09e-02 & 1.00 \\	
			\noalign{\smallskip}\hline
		\end{tabular}	
		\smallskip
		\caption{Temporal accuracy of {\it Case 1} in {\bf Experiment 4}.}
		\label{tab:ex4-u0}
	\end{center}
\end{table}

\begin{table}[H] 
	\small\setlength{\tabcolsep}{2.5pt}
	\begin{center}	
		\begin{tabular}{cccccccccccc}
			\hline\noalign{\smallskip}
			&&  &&\multicolumn{2}{c}{$\alpha = 0.2$}  && \multicolumn{2}{c}{$\alpha=0.4$}
			&& \multicolumn{2}{c}{$\alpha=0.8$}   \\
			\noalign{\smallskip}		\cline{5-6} \cline{8-9}\cline{11-12}\noalign{\smallskip}
			$\sigma$&& $J$ && $\mathcal E_2 $ &{Order}    &
			&  $\mathcal E_2 $ &{Order} &
			&  $\mathcal E_2 $ &{Order} \\
			\noalign{\smallskip}\hline\noalign{\smallskip}
			\multirow{4}{*}{1.5}
			&  & $2^5$  && 3.94e-02 & --   &  &   7.99e-02 & --&  &   1.75e-01& --  \\
			&  & $2^6$  && 2.67e-02& 0.56 &  & 5.81e-02 & 0.46&  &   1.16e-01 & 0.60\\
			&  & $2^7$ && 1.78e-02 & 0.58 &  &  4.02e-02 & 0.53 &  &   7.33e-02& 0.66 \\
			&  & $2^8$ &&  1.16e-02 & 0.62 &  &  2.64e-02 & 0.61 &  &   4.49e-02 & 0.71\\				
			\noalign{\smallskip}\hline\noalign{\smallskip}
			\multirow{4}{*}{2.5}
			&  & $2^5$ && 1.16e-02& --   &  &2.56e-02 & -- &  &   4.90e-02 & -- \\
			&  & $2^6$ & & 5.93e-03 & 0.97 &  &1.31e-02 & 0.97 &  &   2.46e-02 & 0.99 \\
			&  & $2^7$ && 2.94e-03 & 1.01 &  &  6.62e-03 & 0.98 &  &   1.22e-02 & 1.01 \\
			&  & $2^8$&&  1.46e-03 & 1.01 &  &  3.26e-03& 1.02 &  &   6.01e-03 & 1.02 \\
			\noalign{\smallskip}\hline
		\end{tabular}
		\smallskip
		\caption{Temporal accuracy of {\it Case 2} in {\bf Experiment 4}.}
		\label{tab:ex4-f-L2-1}
	\end{center}
\end{table} 

However, for the last case, it seems hard (or even impossible) to obtain growth estimate of the solution, and the accuracy $\mathcal E_2=O(\tau^{1/2})$ can not be improved; see Table~\ref{tab:ex4-f-L2}.

\begin{table}[H] 
	\small\setlength{\tabcolsep}{2.5pt}
	\begin{center}		
		\begin{tabular}{ccccccccccccccc}
			\hline\noalign{\smallskip}
			&&  &&\multicolumn{2}{c}{$\sigma= 1.5$}  && \multicolumn{2}{c}{$\sigma=2.5$}
			&& \multicolumn{2}{c}{$\sigma=5$}
			&& \multicolumn{2}{c}{$\sigma=10$}    \\
			\noalign{\smallskip}	\cline{5-6} \cline{8-9}\cline{11-12}\cline{14-15}\noalign{\smallskip}
			$\alpha$&& $J$ && $\mathcal E_2 $ &{Order}    &
			&  $\mathcal E_2 $ &{Order} &
			&  $\mathcal E_2 $ &{Order} &
			&  $\mathcal E_2 $ &{Order} \\
			\noalign{\smallskip}\hline\noalign{\smallskip}
			\multirow{4}{*}{0.1}
			&&$2^8$  && 2.44e-02 & --   &  &   3.03e-02 & --&  &   4.06e-02& --  &&5.45e-02&--\\
			&&$2^{9}$  && 1.80e-02& 0.43&  & 2.24e-02 & 0.43&  &   3.01e-02 & 0.43&&4.04e-02&0.43\\
			&&$2^{10}$ && 1.32e-02 & 0.45 &  &  1.64e-02 & 0.45 &  &   2.21e-02& 0.45 &&2.96e-02&0.45\\
			&&$2^{11}$ &&  9.56e-03 & 0.47 &  &  1.19e-02 & 0.47 &  &   1.60e-02 & 0.47&&2.14e-02&0.47\\	
			\noalign{\smallskip}\hline\noalign{\smallskip}
			\multirow{4}{*}{0.2}
			&&$2^{9}$  &&2.58e-02& --&  & 3.08e-02 & --&  &   3.92e-02 & --&&4.98e-02&--\\
			&&$2^{10}$ && 1.99e-02 & 0.37 &  &  2.39e-02 & 0.37 &  &   3.04e-02& 0.37 &&3.87e-02&0.36\\
			&&$2^{11}$ &&  1.51e-02 & 0.40 &  &  1.82e-02 & 0.39 &  &   2.32e-02 & 0.39&&2.96e-02&0.39\\	
			&  & $2^{12}$&&  1.10e-02 & 0.46 &  &  1.33e-02& 0.45&  &  1.70e-02 & 0.45 &&2.18e-02&0.44\\
			\noalign{\smallskip}\hline
		\end{tabular}
		\smallskip
		\caption{Temporal accuracy of {\it Case 3} in {\bf Experiment 4}.}
		\label{tab:ex4-f-L2}
	\end{center}
\end{table} 
	
	\appendix 
	
	\section{Some properties of fractional calculus operators} \label{secA1}
	
	\begin{lemma}[\cite{Samko1993}] 	\label{lem:basic-frac} 
		If $ 0 <\alpha, \beta < \infty $, then
		\[
		\I_{0+}^\beta \I_{0+}^\alpha v = \I_{0+}^{\beta +\alpha}v, \quad
		\I_{1-}^\beta \I_{1-}^\alpha v = \I_{1-}^{\beta +\alpha}v,
		\]
		for all $v\in L^1(0,1)$,
		and if $ 0 < \alpha < \beta < \infty $, then
		\[
		\D_{0+}^\beta \I_{0+}^\alpha v = \D_{0+}^{\beta-\alpha}v, \quad
		\D_{1-}^\beta \I_{1-}^\alpha v = \D_{1-}^{\beta-\alpha}v,
		\]
		for all $v\in L^1(0,1)$.
		Moreover, for all $ v, w \in L^2(0,1) $,
		\[
		\dual{\I_{0+}^\beta v,w}_{(0,1)} =
		\dual{v, \I_{1-}^\beta w}_{(0,1)}.
		\]
	\end{lemma}

	\begin{lemma}[\cite{Ervin2006}]	\label{lem:coer} 
		If $ 0 < \gamma < 1/2 $ and $ v,w \in H^{\gamma}(0,T) $, then
		\begin{align*}
			&\dual{\D_{0+}^\gamma v, \D_{T-}^\gamma v}_{(0,T)}
			=\cos\gamma\pi \snm{v}_{H^{\gamma}(0,T)}^2,\\		
			&
			\dual{\D_{0+}^\gamma v, \D_{T-}^\gamma w}_{(0,T)} = \dual{\D_{0+}^{2\gamma} v, w}_{H^\gamma(0,T)} =
			\dual{\D_{T-}^{2\gamma} w, v}_{H^\gamma(0,T)},\\
			&  \cos\gamma\pi \nm{\I_{0+}^\gamma v}_{L^2(0,T)}^2 \leqslant
			\dual{\I_{0+}^\gamma v, \I_{T-}^\gamma v}_{(0,T)}
			\leqslant \sec\gamma\pi\nm{\I_{0+}^\gamma v}_{L^2(0,T)}^2,\\
			&
			\cos\gamma\pi \nm{\D_{0+}^\gamma v}_{L^2(0,T)}^2\leqslant  \dual{\D_{0+}^\gamma v, \D_{T-}^\gamma v}_{(0,T)}
			\leqslant \sec\gamma\pi \nm{\D_{0+}^\gamma v}_{L^2(0,T)}^2.
		\end{align*}
	\end{lemma}

	\begin{lemma}[\cite{Luo2019}] 	\label{lem:interp} 
		If $ v \in {}_0H^\beta(0,1;\dot H^r(\Omega)) \cap {}_0H^\gamma(0,1;\dot H^s(\Omega)) $ with
		$ \gamma,\beta\geqslant 0$ and $ s,r\in\mathbb R$, then for all $ 0 < \theta < 1 $,
		\[
		\begin{aligned}
			{}&			\nm{v}_{
				{}_0H^{\theta \beta + (1-\theta) \gamma}
				(0,1;\dot H^{\theta r + (1-\theta)s}(\Omega))
			} \\
			\leqslant{} &
			C_{\beta,\gamma,\theta} \left(
			\nm{v}_{{}_0H^\beta(0,1;\dot H^r(\Omega))} +
			\nm{v}_{{}_0H^\gamma(0,1;\dot H^s(\Omega))}
			\right).
		\end{aligned}
		\]
		Similarly, if $ v \in {}^0\!H^\beta(0,1;\dot H^r(\Omega)) \cap {}^0\!H^\gamma(0,1;\dot H^s(\Omega)) $ with
		$ \gamma,\beta\geqslant 0$ and $ s,r\in\mathbb R$, then for all $ 0 < \theta < 1 $,
		\[
		\begin{aligned}
			{}&		\nm{v}_{
				{}^0\!H^{\theta \beta + (1-\theta) \gamma}
				(0,1;\dot H^{\theta r + (1-\theta)s}(\Omega))
			} \\
			\leqslant{} &
			C_{\beta,\gamma,\theta} \left(
			\nm{v}_{{}^0\!H^\beta(0,1;\dot H^r(\Omega))} +
			\nm{v}_{{}^0\!H^\gamma(0,1;\dot H^s(\Omega))}
			\right).
		\end{aligned}
		\]
	\end{lemma}

	\begin{lemma}[\cite{Luo2019}] 	\label{lem:coer1} 
		If $\beta\geqslant \gamma>0$, then
		\begin{align*}
			\nm{\D_{T-}^{\gamma} v}_{{}^0\!H^{\beta-\gamma}(0,T)}  \leqslant{}&
			C_{1}\nm{v}_{{}^0\!H^\beta(0,T)}\quad
			\forall\,v \in {}^0\!H^\beta(0,T),\\
			\nm{\D_{0+}^{\gamma} v}_{{}_0H^{\beta-\gamma}(0,T)}  \leqslant{}&
			C_{2} \nm{v}_{{}_0H^\beta(0,T)} \quad
			\forall\,v \in {}_0H^\beta(0,T),
		\end{align*}
		where $C_1$ and $C_2$ depend only on $\gamma$ and $\beta$.
	\end{lemma}

	\begin{lemma}[\cite{Luo2019}] 	\label{lem:regu} 
		If $\beta,\gamma\geqslant 0$, then
		\begin{align*}
			C_{1}\nm{v}_{{}^0\!H^\beta(0,T)}
			\leqslant
			\nm{\I_{T-}^{\gamma} v}_{{}^0\!H^{\beta+\gamma}(0,T)}  \leqslant{}&
			C_{2}\nm{v}_{{}^0\!H^\beta(0,T)}
			\quad
			\forall\,v \in {}_0H^\beta(0,T),\\
			C_{3}\nm{v}_{{}_0H^\beta(0,T)}
			\leqslant
			\nm{\I_{0+}^{\gamma} v}_{{}_0H^{\beta+\gamma}(0,T)}  \leqslant{}&
			C_{4}\nm{v}_{{}_0H^\beta(0,T)}
			\quad
			\forall\,v \in {}_0H^\beta(0,T).
		\end{align*}
		where $C_1,\,C_2,\,C_3$ and $C_4$ depend only on $\gamma$ and $\beta$.
	\end{lemma}
	
	\begin{lemma}[\cite{Luo2019}]
		\label{lem:Nirenberg}
		If $0< \gamma<1/2$, then for all $v \in {}_0H^1(0,1)$,
		\begin{equation*}
			\nm{v}_{C[0,1]} \leqslant C_{\gamma}
			\nm{v}_{{}_0H^1(0,1)}^{(1/2-\gamma)/(1-\gamma)}
			\nm{v}_{{}_0H^\gamma(0,1)}^{1/(2-2\gamma)}.
		\end{equation*}
		Moreover, if $v\in {}_0H^{\gamma}(0,1)$ with $1/2<\gamma\leqslant 1$, then
		for all $ 0<\epsilon\leqslant1 $,
		\[
		\nm{v}_{C[0,1]}\leqslant
		\frac{C_{\gamma}}{\sqrt\epsilon}
		\nm{v}_{{}_0H^{1/2}(0,1)}^{1-\epsilon}
		\nm{v}_{{}_0H^{\gamma}(0,1)}^{\epsilon}.
		\]
	\end{lemma}

	\bibliographystyle{abbrv}
	
\end{document}